\documentclass[12pt]{amsart}
\usepackage{amsmath,amsthm,amsfonts,amssymb}
\usepackage[all]{xy}
\usepackage{amssymb, amsthm, amsmath}
\usepackage[english]{babel}
\usepackage{amsfonts}
\usepackage{tikz}
\usetikzlibrary{matrix}
\DeclareMathOperator{\id}{id}
\DeclareMathOperator{\Hom}{Hom}

\newcommand{\ra}{\rightarrow}

\newcommand{\Z}{\mathbb Z}
\newcommand{\ot}{\otimes}
\newcommand{\mtc}{\mathcal}

\newcommand{\al}{\alpha}
\newcommand{\eps}{\epsilon}
\newcommand{\bn}{\begin}

\newcommand{\sub}{\subsection}

\newcommand{\D}{\Delta}

\numberwithin{equation}{section}
\newtheorem{lemma}[equation]{Lemma}

\newtheorem{defn}[equation]{Definition}
\newtheorem{cor}[equation]{Corollary}
\newtheorem{rem}[equation]{Remark}

\newcommand{\bl}{\begin{lemma}
  }
\newcommand{\nc}{\newcommand}
\nc{\el}{\end{lemma}}

\newcommand{\ch}{\chi}
\newcommand{\mtr}{\mathrm}
\nc{\bwt}{\bowtie}
\newcommand{\ncm}{\newcommand}\newcommand{\gm}{\gamma}
\numberwithin{equation}{section}
\newcommand{\et}{\end{thm}}\newcommand{\bt}{\bn{thm}}
\newcommand{\ep}{\end{prop}}\newcommand{\bp}{\bn{prop}}
\newcommand{\beqarn}{\begin{eqnarray*}}
\newcommand{\eeqarn}{\end{eqnarray*}}
\newcommand{\beqn}{\bn{equation*}}
\newcommand{\eeqn}{\end{equation*}}
\newcommand{\bpf}{\bn{proof}}
\newcommand{\epf}{\end{proof}}
\ncm{\cX}{\mtc{X}}
\ncm{\wt}{\widetilde}
\ncm{\sg}{\sigma}\ncm{\Rep}{\mathrm{Rep}}
\newcommand{\Res}{\mathrm{Res}}
\ncm{\X}{\mathcal{X}}
\ncm{\cA}{\mathcal{A}}
\newcommand{\lb}{\label}

\numberwithin{equation}{section}

\numberwithin{equation}{section}

\ncm{\np}{\newpage}
\ncm{\ebl}{\end{thebibliography}}
\ncm{\bbl}{\begin{thebibliography}}
\ncm{\chd}{_{ _{\ch}}}
\ncm{\ald}{_{ _{\al}}}
\ncm{\cP}{\mathcal{P}}
\ncm{\ei}{e_i}
\ncm{\eij}{e_{i,\;j}}
\ncm{\bne}{\begin{enumerate}}
\ncm{\ene}{\end{enumerate}}\ncm{\bdef}{\begin{defn}}
\ncm{\edf}{\end{defn}}
\ncm{\stab}{\mtr{Stab}}
\ncm{\bc}{\begin{cor}}

\ncm{\ec}{\end{cor}}
\ncm{\er}{\end{rem}}
\ncm{\br}{\begin{rem}}

\ncm{\bd}{\begin{document}}
\ncm{\ed}{\end{document}}
\ncm{\beq}{\begin{equation}}

\ncm{\eeq}{\end{equation}}
\ncm{\cm}{\mathcal{M}}
\ncm{\rep}{\mtr{Rep}}
\ncm{\btw}{\bowtie}
\ncm{\cd}{\mtc{D}}
\ncm{\cop}{\mtr{cop}}

\ncm{\bea}{\begin{eqnarray}}
\ncm{\eea}{\end{eqnarray}}
\ncm{\beanon}{\begin{eqnarray*}}
\ncm{\eeanon}{\end{eqnarray*}}\ncm{\ek}{\eps|_K}\ncm{\diez}{\#}

\ncm{\cC}{\mtc{C}}

\ncm{\cc}{\mtc{C}}

\ncm{\HKer}{\mtr{HKer}}
\ncm{\LKER}{\mtr{LKER}}
\ncm{\aad}{\mtr{ad}}
\ncm{\Dr}{\mtr{D}}
\ncm{\cD}{\mathcal{D}}
\ncm{\G}{\mathcal{G}}
\ncm{\Dc}{\mtc{D}}
\ncm{\E}{\mtc{E}}
\ncm{\fp}{\mtr{FP}}
\ncm{\Vc}{\mtr{Vec}}

\ncm{\cK}{\mtc{K}}
\ncm{\cM}{\mtc{M}}
\ncm{\cE}{\mtc{E}}
\ncm{\cS}{\mtc{S}}
\ncm{\cs}{\cS}
\ncm{\End}{\mtr{End}}
\ncm{\hsa}{Hopf subalgebra of }
\ncm{\ses}{semisimple}
\ncm{\x}{$}
\ncm{\mi}{\mtr{I}}
\ncm{\cZ}{\mtc{Z}}\ncm{\xra}{\xrightarrow}
\ncm{\cb}{\mtc{B}}\ncm{\ca}{\mtc{A}}
\ncm{\irr}{\Irr}\ncm{\Irr}{\mathrm{Irr}}
\ncm{\co}{\mtc{O}}%\ncm{\co}{\mtc{O}}
\ncm{\cg}{{\mtr{K}_0}}\ncm{\ci}{\mtc{I}}
\ncm{\blue}{\textcolor[rgb]{.00, .00, 1.00}}
\ncm{\bb}{\blue}
\ncm{\red}{\textcolor[rgb]{1.00, .00, .00}}
\ncm{\md}{\medbreak}
\ncm{\green}{\textcolor[rgb]{.00, 1.00, .00}}
\ncm{\Gm}{\Gamma}\ncm{\ind}{\mtr{Ind}}\ncm{\res}{\mtr{Res}}
\numberwithin{equation}{section}
\ncm{\bq}{\beq}\ncm{\mto}{\mapsto}\ncm{\opl}{\oplus}
\ncm{\eq}{\eeq}\newcommand{\R}{{\mathcal R}}
\newcommand{\C}{{\mathcal C}}
\ncm{\Ind}{\mtr{Ind}}\ncm{\cz}{\mtc{Z}}\ncm{\ce}{\mtc{E}}\ncm{\ro}{T}
\ncm{\inv}{^{{-1}}}
\title[Fusion categories]
{G - functors arising from categorical group actions on abelian categories}%double crossed products via Clifford theory}
\author{Sebastian  Burciu}
\address{Inst.\ of Math.\ ``Simion Stoilow" of the Romanian Academy
P.O. Box 1-764, RO-014700, Bucharest, Romania}%\address{and}\address{University of Bucharest, Faculty of Mathematics and Computer Science, Algebra and Number Theory Research Center, 14 Academiei St., Bucharest, Romania }
\email{sebastian.burciu@imar.ro} 
\thanks{This work was supported by a grant of the Romanian National Authority for Scientific Research, CNCS-UEFISCDI, project number PN-II-RU-TE-2012-3-0168.}%This work was partially supported by a grant of the Romanian National Authority for Scientific Research, CNCS-UEFSCDI, project number PN II - ID - PCE-2011-3-0168}
\begin{document}
\subjclass[2000]{Primary 16W30, 18D10}

%\keywords{Double crossed products; Drinfeld doubles; Fusion categories; Clifford theory}
\keywords{Group actions on categories; Fusion categories; Grothendieck rings;}
\maketitle
%\contents

\begin{abstract}
 A Mackey type decomposition for group actions on abelian categories is described. This allows us to define new Mackey functors which associates to any subgroup the $K$-theory of the corresponding equivariantized abelian category. In the case of an action by tensor autoequivalences the Mackey functor at the level of Grothendieck rings has a Green functor structure. As an application we give a description of the Grothendieck rings of equivariantized fusion categories under group actions by tensor autoequivalences on graded fusion categories.%Simple objects of an equivariantization under a homogenous group actions are described. A new formula for the tensor product of any two simple objects of an equivariantized fusion category is also given. 
\end{abstract}

%\tableofcontents\newpage
%\input{introd}
\section{Introduction and Main Results}
A Mackey functor or (a $G$-functor) is a family $\{a(K)\}_{K\leq G}$ of abelian groups equipped
with three types of maps: induction, conjugation, and restriction, satisfying some certain compatibility axioms, see for example \cite{gr71}. Typical examples, include among others, the cohomology groups $\{H^{n}(K, M)\}_{K\leq G}$ and the character rings $\{R(K)\}_{K \leq G}$ are both G-functors.
\md
It is shown in \cite{35, 50, 8} that the class group of the ring of integers of the fixed field $\{k^{H}\}_{H \leq G}$ where $G$ is a group of automorphisms of a number field $k$ is a Mackey functor.
In a somehow different direction, these results were extended in \cite{kevin} by showing that $\{K_{i}(S^{H })\}_{H\leq G}$ is a $G$-functor, whenever $R\subseteq S$ is a Galois extension of commutative rings with Galois group $G$.
\md
The main goal of this paper is to construct in the same spirit new $G$-functors arising from categorical group actions on categories. More precisely, if a finite group $G$ acts on an abelian category $\cc$ we show that the $K$-theory $\{K_{i}(\cc^{H })\}_{H\leq G}$ is also $G$-functor %for all $i \geq 0$ 
(see Theorem \ref{main2}). If $\cc$ is a tensor category %(over a field $k$) 
and $G$ acts by tensor autoequivalences then the $\{K_{0}(\cc^{H })\}_{H\leq G}$ is a Green functor.
\md
Let $G$ be a finite group acting on the abelian category $\cc$.
For any subgroup $H$ of $G$ the left adjoint functor of the forgetful functor $\res^{G}_{H}:\cc^{G}\ra  \cc^{H}$ was recently described in \cite{buna}.  This functor is denoted by $\Ind^{G}_{H}:\cc^{H}\ra  \cc^{G}$ and can be regarded as a generalization of the induction functor from $\rep(H)$ to $\rep(G)$. 
\md
%In particular we prove a Mackey type  decomposition in the  general setting obtained from an action of a finite group $G$ on an arbitrary $k$-linear category $G$,  see Theorem \ref{macky} below.
%\md The classical Mackey decomposition for induced representations of a finite group can be stated as a compatibility condition between the induction and restriction functors for the categories of representations of various subgroups of $G$. This decomposition plays an important role in the theory of representations of finite groups, especially for the modular case \cite{benson}.
%\blue{It is easy to notice that the same formula for $\Ind^{G}_{H}$ works for a group action on any tensor category $\cc$.}
\md
Our first main result is the following: 
 \bt\lb{main2}
Let $G$ be a finite group acting on the abelian category $\cc$. Then for all $i \geq 0$ the association $H \mapsto K_{i}(\cc^{H}) $ defines a Mackey functor $M_{i}$ with the following structure maps:
\bne
\item Restriction
$
R^{H}_{K} :K_{i}(\cc^{H}) \ra K_{i}(\cc^{K})
$
is the map induced by the forgetful functor $\res^{H}_{K}:
\cc^{H}\ra \cc^{K}$ and 
\item Induction $
I^{H}_{K} :K_{i}(\cc^{H}) \ra K_{i}(\cc^{K})
$
is the map induced by the induction functor $\ind^{H}_{K}:\cc^{K}\ra \cc^{H}$,
\item Conjugation $
c_{H, x} :K_{i}(\cc^{H})\ra K_{i}(\cc^{\;^{x}H})
$
is the map induced by the functor $T^{x}:\cc^{H}\ra \cc^{\;^{x}H}$.
\ene
\et
If $\cc$ is a tensor category  over a fixed field $k$ and the action of $G$ on $\cc$ is by tensor autoequivalences then we show that for $i=0$ the above Mackey functor $M_{0}$ is in fact a Green functor.
\bt\label{grf}
Let $\cc$ be a tensor category over a fixed field $k$ and $G$  a finite group acting  on $\cc$ by tensor autoequivalences. Then $H \mapsto K_{0}(\cc^{H})$ defines a Green functor on $G$ over $k$.
\et
The proof of the above results uses a Mackey type decomposition for the above induced functor when restricted to various subgroups:
\bt\lb{macky}
Suppose that a $G$ acts on an abelian category 
$\cc$ via $T:\underline{G} \to \underline{ \mtr{Aut}}(\cc)$. Let $K$ and $H$ be any two subgroups of $G$ and $M\in \cc^H$. Then
\beq
\res_K^G(\ind^G_H(M))\simeq\bigoplus_{x \in D}\ind_{K \cap \;^xH}^K(\res^{\;^xH}_{\;^xH\cap K}(T^x(M)))
\eeq
where $D$ is a complete set of representatives for the space of double cosets $K\backslash G \slash H$ and $\;^{x}H:=xHx^{{-1}}$.
\et
%For a finite group $G$ denote by $\cs(G)$ denotes the lattice of all subgroups of $G$.\md
% Our second main result is the above mentioned  construction of Green functors associated to any action of a finite group on a tensor category. 
 %\blue{put them as identities on functors and they induce identities on K-theory,}
%As a Corollary of this theorem we obtain that the map  $\cs(G)\ra k-\text{alg}$ given by $H\mapsto K_{0}(\cc^{H} )$ is a Mackey functor where. This becomes a Green functor if the action of $G$ is  by tensor automorphisms. 
%\blue{As an application of this theory we construct vertices and sources for equivariant objects under a group action on tensor categories. In particular we obtain a description of indecomposable objects of an equivariantization. Quantum groups $u_{q}(g)$, Clifford theory on general satisfying property $P$.}\md
%We prove a version of Frobenius reciprocity  for these pairs of adjoint functors.\blue{this by definition}
\md
The last section of this paper is concerned with the structure of the Grothendieck rings of equivariantized fusion categories under certain group actions. These actions are called in the paper {\it coherent actions} since they respect a certain natural compatibility with  the grading of the fusion category.  We show that under a coherent action the Grothendieck group of an equivariantization has the structure of the rings introduced in \cite{scoh}. These rings are also considered in \cite{bouc} as Green rings obtained from Dress construction from other given Green rings. Examples include, among others, the crossed Burnside rings, the Hochschild cohomology rings of crossed products, and the Grothendieck rings of (twisted) Drinfeld double of finite groups. Using Theorem \ref{grf} it is shown in this paper that the Grothendieck ring of  the Drinfeld center of any fusion category has this structure.
\md

This paper is organized as follows. Section \ref{prelim} recalls some basic results on abelian categories and group actions on them. The construction of the adjoint functor $\Ind^{G}_{H}$ mentioned above is recalled in this section. In Section \ref{tim} the proof of the main Theorem \ref{macky} is presented. Section \ref{green} is devoted to the proof of Theorem \ref{main2} and Theorem \ref{grf}. In this section we recall the definition of $G$-functors and Green functors. Last section is devoted to the study of the Grothendieck rings of equivariantizations of graded fusion categories. Coherent group actions on graded fusion categories are introduced in this section. In Proposition \ref{desgnn} we give a new description for the simple objects of an equivariantization under a coherent group action, generalizing results of \cite[Proposition 2.7]{gnn}. In this section we also recall the ring structures introduced in \cite{scoh} and prove that the Grothendieck groups of equivariantizations under coherent actions has this type of  structure.% \bb{In section \ref{hmg} the description of simple objects of a homogenous equivariantization are described. A tensor product formula is also given. In section \ref{withsp} it is shown that the Grothendieck ring of a homogenous equivariantization has the structure of a Green ring as in \cite{scoh}.}

%We work on an arbitrary field $k$.
\section{Group actions on categories}\lb{prelim}
%\subsubsection{Abelian categories}
%Recall that an additive category is called preabelian if every morphism has both a kernel and a cokernel. Moreover, a preabelian category is called abelian if every monomorphism and every epimorphism is normal. This means that every monomorphism is a kernel of some morphism, and every epimorphism is a cokernel of some morphism.  %All the abelian categories through this paper are supposed to admit finite direct sums of objects.
%\md
%}Recall that a functor $F:\cc \ra \cd$ between two abelian categories is exact if and only if it preserves direct sums and kernels.
\subsection{Tensor categories} 
\subsubsection{$k$-linear categories}
Fix a commutative ring $k$. Recall that a $k$-linear category is an abelian category in which the hom-sets are $k$-vector spaces, the compositions are $k$-bilinear. A $k$-linear functor between $k$-linear categories is a functor which is linear on all hom-spaces.
\md
Recall that an essentially small $k$-linear category is said to be locally
finite \cite{EGNO} if for any two objects $X ,Y \in \cc$, the space $\Hom_{\cc}(X, Y )$ is finite dimensional and every  object in $\cc$ has finite length.
%(ii) every.
\subsubsection{Tensor categories} 
Let $\cc$ be a $k$-linear rigid
monoidal category. Then $\cc$  is called a tensor category over $k$ if
the tensor bifunctor  is bilinear on morphisms and $\End(1)=k$ (see \cite{ENO}).
\subsubsection{Tensor functors and natural tensor transformations}
%then we will call $\cc$ a tensor category \subsection{Group actions on $k$-linear categories} Let us recall the definition of the group action on a $k$-linear category from \cite{tamb}. See also \cite{gait, N} for the tensor settings.
 %finite tensor category \cite{gait, N, tamb}.
 \md
Recall that a {\it unitary tensor functor} $F:\cc \ra \cd$ between two tensor categories is a $k$-linear functor $F$ together with a natural transformation $F_{2}:F(- \; \ot \; -)\ra F(-)\ot F(-)$ satisfying several compatibility axioms (see for example \cite{ENO}). 

In particular the naturally of $F_{2}$ with respect to the morphisms can be written as
\beq\label{natf}
(F(u)\ot F(v))F_{2}^{M, N}=F_{2}^{M', N'}F(u \ot v)
\eeq
for all morphisms $M\xra{u} M'$ and $N\xra{v} N'$ in $\cc$.
Composition of two tensor functors $\cc \xra{G} \cd \xra{F}\ce$ is also a tensor functor with 
\beq\label{comp}
(F\circ G)_{2}^{M, N}:=F_{2}^{G(M), G(N)}\circ F(G_{2}^{M, N})
\eeq
A {\it natural tensor transformation }$\tau: F\ra G$ between two tensor functors is a natural transformation satisfying the following compatibility condition:
\beq\label{tenant}
G^{M, N}_{2}\tau_{M\ot N}=(\tau_{M}\ot \tau_{N})F^{M, N}_{2}
\eeq
for any objects $M, N \in \cc$. 
\subsubsection{Fusion categories}
Let $k$ be an algebraically closed field. A fusion category over $k$ is a rigid semisimple $k$-linear tensor category $\cc$ with finitely many simple
objects and finite dimensional spaces of morphisms such that the unit object of $\cc$ is simple.
\subsection{Group actions on abelian categories}%$k$-linear categories}
Let $\C$ be an abelian category. Denote by $ \underline{ \mtr{Aut}}(\cc)$ the category whose objects are {\it exact} autoequivalences of $\cc$ and morphisms are natural transformations between them. Then $ \underline{ \mtr{Aut}}(\cc)$ is a monoidal category where the tensor product is defined as the composition of autoequivalences.
\md
For a finite group $G$ let $\mathrm{Cat}(G)$ denote the monoidal category whose objects
are elements of $G$, the only morphisms are the identities, and the tensor product is given by multiplication in $G$.\md
An action of a finite group $G$ on $\cc$ consists of a unitary monoidal functor $T:\mathrm{Cat}(G) \ra \underline{ \mtr{Aut}}(\cc)$. Thus, for every $g \in G$, we have a functor $T^g: \C \to \C$ and a collection of natural isomorphisms 
\beqn T^{g,h}_2 : T^g T^h \to T^{gh}, \quad g, h \in G,\eeqn which give the tensor structure of $G$. The tensor unit of $T$ is denoted by $T_{0}:\id_{{\cc}}\ra T^{1}$ where $1 \in G$ is the unit of the group $G$. \md By the definition of the tensor functor, the tensor structure $T_{2}$ satisfies the following conditions:
\begin{align}\label{ro-2} & (T^{gh, l}_2)_M \, (T^{g,h}_2)_{T^l(M)} =
(T^{g, hl}_2)_M \, T^a((T^{h,l}_2)_M), \\ & (T^{g,
e}_2)_M T^{g}({T_0}_{ _{M}}) = (T^{e, g}_2)_M
(T_0)_{T^{g}(M)},
\end{align}
for all objects $M \in \C$, and for all $g, h, l \in G$. See \cite[Subsection 4.1]{DGNO}. 
Note that by the naturality of $T^{g, h}_2$, $g, h \in G$,  can be written as
\begin{equation}\label{nat-ro}T^{gh}(f) \, (T_2^{g, h})_{N} = (T_2^{g, h})_M \, T^gT^h(f),
\end{equation}
 for every morphism $f: M \to N$ in $\C$.\md

We shall assume in what follows that $T^1 = \id_\C$ and
$T_0$, $T^{g, 1}_2$, $T^{1, g}_2$ are also identities. We say that $G$ acts {\it $k$-linearly} on the $k$-linear category $\cc$ if $T^{g}$ is a $k$-linear autoequivalence for any $g \in G$.
\bn{example}\label{automact}
Suppose that $G$ acts as a ring automorphisms on a $k$-algebra $S$. Then $G$ acts on $S$-mod via  the following action: $T^{g}(M)=M$ as abelian groups and the $S$-action on $T^{g}(M)$ is given by $s.\;^{g}m:=(g^{-1}.s)m$. In this case one can take $(T^{g, h}_{2})_{M}=\id_{M}$ for all $g,h \in G$.
\end{example}
\subsection{On the equivariantized category}
\medbreak Suppose that $G$ acts on the abelian category $\cc$. Let $\C^G$ denote the corresponding
\emph{equivariantized} category. Recall that $\C^G$ is an abelian category whose objects are $G$-equivariant
objects of $\C$. They consist of pairs $(M, \mu)$, where $M$ is an object
of $\C$ and $\mu = (\mu_{M}^g)_{g \in G}$ is a collection of isomorphisms $\mu_{M}^g:T^g(M) \to M$ in $\cc$ satisfying the following:
\begin{equation}\label{deltau} \mu_{M}^g T^g(\mu_{M}^h) = \mu_{M}^{gh} (T^{g,
h}_2)_M, \quad \forall g, h \in G, \quad \quad
\mu_{M}^{1}{T_0}_M=\id_M.\end{equation} \medbreak We say that an object $M$ of $\C$ is
\emph{$G$-equivariant} if there exists such a collection $\mu =
(\mu^g)_{g \in G}$ so that $(M, \mu) \in \C^G$. Note that the equivariant structure $\mu$ is not necessarily unique.
\md 
A morphism $f : (M, \mu_{M})\ra (N, \mu_{N})$ in $\cc^{G}$ is a morphism in $f:M\ra N$ such that 
\beq\label{morf}
\mu^{g}_{N}T^{g}(f)=f\mu^{g}_{M}.
\eeq

\bn{example}\label{automequiv}
It is easy to verify that in the case of the previous example one has that $(S$-mod)$^{G}\simeq S\#kG$-mod, the category of $S\#kG$-modules.
\end{example}
\subsection{Induction functors as left adjoints of restriction functors} Suppose that a finite group $G$ acts on the abelian category $\cc$ and let $H\leq G$ be a subgroup. Let $\R$ be a set of representative elements for the left cosets
of $H$ in $G$. Thus one can write $G$ as a disjoint union $G = \cup_{t \in \R}tH$. 
Set, for all $(V, \mu) \in \cc^{H}$, 
\bq\lb{indx}
\ind_H^G(V, \mu): =(\oplus_{t \in \R}T^t(V), \nu)
\eq 
where for all $g \in G$ the equivariant structure of $\nu^{g}: \bigoplus_{t \in \R}T^gT^t(V) \to \bigoplus_{t \in \R}T^t(V)$ is defined componentwise by the formula 
\begin{equation}\label{muind}
\nu^{g, t} = T^s(\nu^h)
(T_2^{s, h})^{-1} T_2^{g, t}: T^g T^t(V) \to
T^s(V).\end{equation} 
Here the elements $h \in H$ and $s \in \R$ are uniquely determined by the relation $gt = sh$.

Note that the proof of  \cite[Proposition 2.9]{buna} works in any abelian category, therefore  $\Ind^{G}_{H}$ is a left adjoint functor of $\Res^{G}_{H}$.

\subsection{Action by tensor equivalences.}
Suppose that $\C$ is a tensor category over $k$  and consider $\underline {\mtr{Aut}}_{\otimes}(\cc) $ the full subcategory of $\underline {Aut}(\cc)$ consisting of $k$-linear tensor autoequivalences of $\cc$.

Let $T: \underline G \to \underline {\mtr{Aut}}_{\otimes} \C$ be an action of $G$ on $\C$ by \emph{tensor autoequivalences},
that is, $T^{g}$ is a tensor auto equivalence for all $g \in G$. Thus $T^g$ is endowed with a monoidal structure
$(T_2^g)^{M, N}: T^g(M \otimes N) \to T^g(M) \otimes
T^g(N)$, for all $M, N \in \C$ and $T_2^{g, h}: T^g T^h \to T^{gh}$ are natural isomorphisms of tensor functors, for all $g, h \in G$. Thus, for
all $g, h\in G$ and $M, N \in \C$ the following relation holds:
\begin{equation}\label{tensor-rho} {(T_2^{gh}})^{M, N} {(T_2^{g, h})}_{M \otimes N} =  ({(T_2^{g, h})}_{M} \otimes {(T_2^{g, h})}_{N}) \, {(T_2^{g})}^{T^h(M), T^h(N)} \, {T^g((T_2^{h})}_{M, N}).
\end{equation} 
\section{Proof of Theorem \ref{macky}}\lb{tim}
\bp\lb{conj} 
Let $G$ be a finite group acting on the abelian category $\cc$. Let $H$ be any subgroup of $G$ and $x \in G$. If $M=(V, \mu) \in \cc^H$ then $T^x(M):=(T^{x}(V), \;^{x}\mu)\in \cc^{\;^xH}$ with the equivariant structure $\:^{x}\mu^{xhx\inv}_{T^{x}(V)}:T^{xhx^{-1}}(T^x(V))\ra  T^x(V)$ given as follows:
\beq\label{conjmu}
T^{xhx^{-1}}(T^x(V))\xra{(T_2^{xhx^{-1}, \;x})_{V}}T^{xh}(V)\xra{(T_2^{x,h})^{-1}_{V}}T^x(T^h(V))\xra{T^x(\mu^{h}_{V})} T^x(V).
\eeq
\ep
ces
\bpf
It is enough to verify Equation \eqref{deltau} which is equivalent to the diagram made of solid arrows below being commutative.

 Note that compatibility conditions (\ref{ro-2})-(\ref{deltau}) of the action of $G$ imply the commutativity diagram after inserting the dashed arrows. Indeed, the bottom right trapeze $(5)$ is commutative by applying $T^{x}$ to the equivariantized condition (\ref{ro-2}) for $V\in \cc^{H}$. The adjacent trapeze $(6)$ is commutative by the naturallity of  $T_{2}^{{x,h}} $ with respect to the morphism $\mu_{V}^{l}$. The rectangle $(4)$ is commutative due to the associativity of the action, Equation \eqref{ro-2}. The parallelogram $(2)$ is commutative due to the associativity of the action, Equation \eqref{ro-2}. Diagram $(3)$ is commutative due to the  naturallity of the natural transformation $T^{xhx^{-1}, x}_{2}$ with respect to the morphism $T^{l}(V)\xra{\mu_{V}^{l}} V$. Diagram $(1)$ is commutative due to the associativity of the action,  Equation \eqref{ro-2}.
 \begin{center}
{\tiny
\begin{tikzpicture}
  \matrix (m) [matrix of math nodes, row sep=3em,column sep=1.15em,minimum width=3em]
  {
     T^{xhx\inv}(T^{xlx\inv}(T^{x}(V))) &  & T^{xhx\inv}(T^{xl}(V)) &  & T^{xhx\inv}(T^{x}(T^{l}(V))) \\%1
     & (1) & & (2) & \\%2
    T^{xhlx\inv} (T^{x}(V)) & &  & &  T^{xhx\inv}(T^{x}(V))\\%3
     & &  &  (3) & \\%4
     T^{xhl}(V)& & T^{xh}(T^{l}(V)) &  & T^{xh}(V)\\%5 
     & (4) & & (6) &  \\%6
    T^{x}(T^{hl}(V)) & & T^{x}(T^{h}(T^{l}(V))) & & \\%7
    & (5) & & & \\%8
     T^{x}(V) & & & &  T^{x}(T^{h}(V)) \\
     };%9
  \path[-stealth]
    (m-1-1) edge node [left] {$(T_{2}^{xhx\inv, xlx\inv})_{V}$} (m-3-1)%vertical A8
    (m-1-3) edge [dashed] node [right] {$(T_{2}^{xhx\inv, xl})_{V}$} (m-5-1)%verticalA10
    (m-1-1) edge node [above] {$T^{xhx\inv}(T_{2}^{xlx\inv , x})_{V}$}(m-1-3)%A1horiz
  (m-3-5) edge node [right] {$(T_{2}^{xhx\inv, x})_{V}$} (m-5-5)%verticalA5
  (m-1-3) edge node [right] [above]{$T^{xhx\inv}(T_{2}^{x, l})_{V}$} (m-1-5)%horizA2
         (m-1-5) edge node [right] {$T^{xhx\inv}(T^{x}(\mu_{V}^{l}))$} (m-3-5)%horizA3
          (m-1-5) edge [dashed] node [left] {$(T_{2}^{xhx\inv, x})_{T^{l}(V)}$}
          (m-5-3)%horizA4
           (m-5-3) edge [dashed] node [below] {$T^{xh}(\mu^{V}_{l})$}(m-5-5)%horiz A6
            (m-5-3) edge [dashed] node [right] {$(T_{2}^{x,h})_{T^{l}(V)}$}(m-7-3)%A14
            (m-7-1) edge [dashed] node [below] {$T^{x}((T_{2}^{h,l})_{V})$} (m-7-3)%A13
          (m-5-1) edge [dashed] node [below] {($T_{2}^{xh, l})^{-1}_{V}$} (m-5-3)%horizA11
    (m-3-1) edge node [left] {$(T_{2}^{xhlx\inv , x})_{V}$} (m-5-1)%verticalA9
    (m-5-5) edge node [right] {$(T_{2}^{x,h})_{V}$} (m-9-5)%verticalA7
     %verticaql
   (m-5-1) edge node [right] {$(T_{2}^{x, hl})^{-1}_{V}$} (m-7-1)%A12
    %(m-3-1) edge [dashed] node [right] {$\simeq$} (m-2-2)
  (m-7-3) edge [dashed] node [right] {$T^{x}(T^{h}(\mu^{V}_{l}))$} (m-9-5)%A15
   (m-7-1) edge node [right] {$T^{x}(\mu^{hl}_{V})$} (m-9-1)%A16
      (m-9-5) edge node [below] {$T^{x}(\mu^{h}_{V})$} (m-9-1)%A17
         ;
\end{tikzpicture}
}
\end{center}

\epf
\subsection{Proof of Theorem \ref{macky}}
We are now ready to give a proof for Theorem \ref{macky}. 
\bpf  Suppose that $M=(V, \mu_{V}) \in \cc^{H}$. 
Then \bq
 \ind^{G}_{H}(M)=(\bigoplus_{x \in G/H}T^{x}(V), 
 \;\mu_{ \ind^{G}_{H}(M)}).
 \eq
 where the equivariant structure is  given on components 
 \beq
 \mu^{g, x}_{ \ind^{G}_{H}(M)}:T^{g}(T^{x}(V))\xra{} T^{gx}(V)=T^{yh}(V)\xra{} T^{y}T^{h}(V)\xra T^{y}(V)
 \eeq
Since  $G=\sqcup_{x\in D}KxH$ formula (\ref{indx}) for induced objects becomes
 \bq
 \ind^{G}_{H}(M)=(\bigoplus_{x \in D}(\bigoplus_{a \in KxH/H}T^{a}(V)), \mu_{ \ind^{G}_{H}(M)}).
 \eq
For any $x \in D$ let 
 \beq 
 V_x:=  \bigoplus_{a \in KxH/H}T^a(V).
 \eeq
Using formula (\ref{muind}) it can be easily verified that the induced equivariant structure $\nu:=\mu_{\ind^{G}_{H}(M)}$ of $\ind^{G}_{H}(M)$ sends the component $V_{x}$ to itself. Indeed, for any $a \in KxH$ if $l \in K$ and $la=bh$ with $h \in H$ then $b=lah\inv \in KaH=KxH$.
 \md
It follows that  $M_{x}:=(V_{x}, \nu|_{K}) \in \cc^{K}$ and then one can write 
\beq
\res_K^G(\ind^G_H(M))\simeq\bigoplus_{x \in D}M_x%\ind_{K \cap \;^xH}^K(\res^{\;^xH}_{\;^xH\cap K}(T^x(M))).
\eeq
as objects of $\cc^{K}$.

Let $V'_{x}:=\ind_{K \cap \;^xH}^K(\res^{\;^xH}_{\;^xH\cap K}(T^x(M))$ and \\ $M'_{x}:=(V'_{x},\; \mu_{\ind_{K \cap \;^xH}^K(\res^{\;^xH}_{\;^xH\cap K}(T^x(M)))}) \in \cc^{K}$ be the corresponding equivariant object.

Then it is enough to show that  for all $x \in D$ one has
\beq
M_x\simeq M'_{x}
 \eeq
 as objects in $\cc^K$.
\md
Note that there is a bijection between the following sets of left cosets
\beq
K/K\cap \;^xH \ra KxH/H
\eeq
 given by $a(K\cap \;^xH)\mapsto axH$.
 \md
 This enables us to write
 \bq\lb{d2}
 V_{x}=\bigoplus_{a \in K/K\cap \;^{x}H}T^{ax}(V)
 \eq
 \md

On the other hand using formula (\ref{indx}) it follows that 
\bq
M'_{x}\simeq \oplus_{a \in K/K\cap\;^{x}H}T^{a}(T^{x}(V)).
\eq
It will be shown that \beqn M'_{x} \xra{f_{x}:=\bigoplus_{a \in K/K\cap \;^{x}H}(T^{a, x}_{2})_{V} }M_{x}\eeqn is a morphism in $\cc^{K}$. 
In order to do this it remains o check that the morphism
\bq
(T^{a,x}_{2})_{V}:T^{a}(T^{x}(V))\ra T^{ax}(V)
\eq
is compatible with the two equivariant structures of the objects $M_{x}$ and $M'_{x}.$ \md Suppose that $l \in K$ and $a, b \in K /K\cap\;^{x}H$ with $la=bl'$ for some $l' \in K\cap \;^{x}H$. Then $lax=bl'x=bxh$ with $h=(x\inv l'x) \in H$.\md

Using again formula (\ref{indx}) the equivalent structure $\widetilde{\nu}$ of $M'_{x}$ on the component $T^{a}(V)$ from Equation (\ref{d1}) is given on components by
 \beq\lb{d1}
 \widetilde{\nu}^{l, a}_{M'_{x}}:T^{l}(T^{a}(T^{x}(V)))\xra{(T_{2}^{l,a})_{T^{x}(V)}}T^{la}(T^{x}(V))=\eeq \vskip -0,5cm \beqn =T^{bl'}(T^{x}(V))\xra{(T^{2}_{b,l'})^{-1}_{T^{x}(V)}} T^{b}(T^{l'}(T^{x}(V)))\xra{T^{b}((\;^{x}\mu_{T^{x}(V)})^{l'})} T^{b}(T^{x}(V)).
 \eeqn
On the other hand the equivariant structure $\nu^{l,a}$ of $M_{x}$ on the component $T^{ax} (V)$  from Equation (\ref{d2}) has the following formula
 \bq
\nu^{l,a}_{M_{x}}:T^{l}(T^{ax}(V))\xra{(T_{2}^{la,x})_{V}} T^{lax}(V)=\eeq\vskip -0,5cm \beqn =T^{bxh}(V)\xra{(T_{2}^{bx,h})^{-1}_{V}} T^{bx}(T^{h}(V))\xra{T^{bx}(\mu^{V}_{h})} T^{bx}(V).
 \eeqn
Using the compatibility properties for the action of $G$ it is easy to verify that these two equivariant structures coincide under $(T_{2}^{a,x})_{V}$, i.e. the following diagram commutes.
\begin{center}
{\Tiny
\begin{tikzpicture}
  \matrix (m) [matrix of math nodes, row sep=5em,column sep=2em,minimum width=4em]
  {T^{l}(T^{a}(T^{x}(V))) & T^{la}(T^{x}(V))=T^{bl'}(T^{x}(V)) & & T^{b}(T^{l'}(T^{x}(V)))  \\%1
  T^{l}(T^{ax}(V))  &  &  & & \\ %2
   T^{lax}(V)=T^{bxh}(V)=T^{bl'x}(V) & & & T^{b}(T^{l'x}(V))=T^{b}(T^{xh}(V))\\%3
  T^{bx}(T^{h}(V)) &  & & T^{b}(T^{x}(T^{h}(V)))\\%4
  T^{bx}(V) &  & &T^{b}(T^{x}(V)) \\%5
  \\%6
};
  \path[-stealth] 
  (m-1-1) edge node [left] {$T^{l}(T_{2}^{a,x})_{V}$} (m-2-1) 
  %(m-1-2) edge [dashed] node [right] {$(T_{2}^{la,x})_{V}$} (m-3-4) 
%  (m-1-4) edge node [right] {$(T_{2}^{b,x})_{V}$} (m-3-4)
  (m-3-4) edge node [right] {$T^{b}(T_{2}^{x,h})_{V}$} (m-4-4)
   (m-4-4) edge node [right] {$T^{b}(T^{x}(\mu^{V}_{h}))$} (m-5-4)
      (m-5-4) edge node [above] {$(T_{2}^{b, x})^{-1}_{V}$} (m-5-1)
  (m-2-1) edge node [left] {$T^{b}(\mu^{V})_{V}$} (m-3-1)
  (m-1-4) edge  node [left] {$T^{b}(T^{l', x}_{2})_{V}$} (m-3-4)
   (m-1-2) edge [dashed] node [right]{$(T^{la,x}_{{2}})_{V}=(T^{bl',x}_{{2}})_{V}$} (m-3-1) 
   (m-4-1) edge node [left]{$T^{bx}(\mu^{V}_{h})$} (m-5-1)

    (m-1-1) edge node [above]{$(T^{l,a}_{{2}})_{T^{x}(V)}$} (m-1-2) 
     (m-1-2) edge node [above]{$(T^{b, l'}_{{2}})_{T^{x}(V)}$} (m-1-4) 
    % (m-1-3) edge node [above]{$T^{b}(T_{2}^{{l', \;x}})_{V}$} (m-1-4) 
      (m-3-1) edge [dashed]  node [above]{$(T_{2}^{b, xh})_{V}$} (m-3-4) 
           (m-3-1) edge  node [left]{$(T_{2}^{bx, h})^{-1}_{V}$} (m-4-1) 
           (m-4-1) edge [dashed]  node [above]{$(T_{2}^{b,x})^{-1}_{T^{h}(V)}$} (m-4-4) 
    ;
\end{tikzpicture}
}
\end{center}
The bottom rectangle is commutative by the naturally of $T^{b, x}_{2}$ with respect to the morphsims, Equation \eqref{nat-ro}. The above rectangle is commutative due to Equation \eqref{ro-2}, the associativity of the action. The upper left diagram is commutative by the same reason. The upper right trapeze is commutative by associativity of the action,  Equation \eqref{ro-2}.
\epf
\br Note that $\Rep(G)$ can be regarded as the equivariantization $\Vc^{G}$ of the trivial action of $G$ on $\cc=\Vc$, \cite{ENO}. In this case the previous theorem recovers the usual Mackey decomposition for representations of finite groups.
\er

\section{$G$ - functors associated to equivariantizations}\lb{green}
This section is devoted to the the proof of Theorem \ref{main2}.

\bl\lb{commir} Let $G$ be a finite group acting on the abelian category $\cc$. Suppose that $A$ and  $B$ are two subgroups of $G$. Then the following identities hold for any $x \in G$:
\beq
T^{x}\circ \res_{B}^A=\res_{\;^xB}^{\;^xA} \circ T^{x}
\eq	
as functors from $\cc^{B}$ to $\cc^{\;^{x}A}$. Also
\bq
T^{x}\circ \ind_{B}^A\simeq \ind_{\;^xB}^{\;^xA} \circ T^{x}
\eq 
as functors from $\cc^{B}$ to $\cc^{\;^{x}A}$
\el
\bpf
The first identity is straightforward. We will verify the second identity.
\md
Let $\mtc{R}$ be a set of representative elements for the left cosets $A/B$.
Suppose that $(M, \{\mu^{b}_{M}\}_{b\in B})\in \cc^B$. Then using formula (\ref{muind}) it follows that
\bq
\ind_B^A(M)=(\opl_{r \in \mtc{R}}T^r(M), \nu_{\ind_B^A(M)})
\eq
where $\nu^{a}$ is defined on the components as follows
\bq
\nu^{a, r}:T^a(T^r(M))\xra{T^{a,r}_{2}(M)} T^{ar}(M)\xra{{(T_{2}^{r', b}})_{M}^{-1}} T^{r'}(T^b(M))\xra{T^{r'}(\mu_M^{b})
} T^{r'}(M)\eq
where $r' \in \mtc{R}$ and $b \in B$ are determined by $ar=r'b$. On the other hand, using Proposition \ref{conj} as an object of $\cc^{\;^{x}A}$ one has that
\bq
T^{x}(\ind^{A}_{B}(M))=(\oplus_{r \in A/B}T^{x}(T^{r}(M)), \;^{x}\nu_{T^{x}(\ind^{A}_{B}(M))})
\eq
with the equivariant structure $\;^{x}\nu$ given on the components  by:
\beqn
 \;^{x}\nu^{xax^{-1}, r}_{M}:T^{xax\inv}(T^{x}(T^{r}(M))) \xra{({T_{2}^{{xax\inv, x}} })_{T^{r}(M)}} T^{xa}(T^{r}(M))\xra{{(T_{2}^{x,a})^{-1}}_{T^{r}(M)}}\eeqn \beqn  \xra{{(T_{2}^{x,a})^{-1}}_{T^{r}(M)}}  T^{x}(T^{a}(T^{r}(M))) \xra{T^{x}(T_{2}^{a, r})^{-1}_{M}} T^{x}(T^{ar}(M))\xra{T^{x}(T_{2}^{r', b})^{{-1}}_{M}}\eeqn \beqn \xra{T^{x}(T_{2}^{r', b})^{{-1}}_{M}} T^{x}(T^{r'}(T^{b}(M)))\xra{T^{x}(T^{r'}(\mu^{b}_{M}))} T^{x}(T^{{r'}}(M))
\eeqn

On the other hand note that 
\bq
\ind_{\;^xB}^{\;^xA}T^{x}(M)=(\oplus_{r \in \mtc{R}}T^{xrx\inv}(T^{x}(M)), \;\eta^{}_{\ind_{\;^xB}^{\;^xA}T^{x}(M)})
\eq
since $x\mtc{R}x\inv$ is a set of representative for the left cosets of $\;^{x}A/\;^{x}B$. Using again formula (\ref{muind}) it follows that the equivariant structure of $\ind_{\;^xB}^{\;^xA}T^{x}(M)$ on the components is given as follows:

\beqn
\eta^{xax^{-1}, \;xrx^{-1}}_{M} :T^{xax\inv}(T^{xrx\inv}(T^{x}(M)) \xra{({T_{2}^{{xax\inv, xrx\inv}} })_{T^{x}(M)}} \eeqn \beqn T^{xarx\inv}(T^{x}(M)) \xra{{(T_{2}^{xr'x\inv,xbx\inv})^{-1}}_{T^{x}(M)}}  T^{xr'x\inv}(T^{xbx\inv }(T^{x}(M)))\eeqn \beqn \xra{T^{xr'x\inv}(\;^{x}\mu_{M})_{xbx\inv}} T^{xr'x\inv}(T^{x}(M))
\eeqn
Define the natural transformation \beqn F:=\oplus_{r \in \mtc{R}}(T_{2}^{xrx\inv, x})^{-1}(T_{2}^{x,r}):T^{x}(\ind_{B}^{A})\ra \ind_{\;^xB}^{\;^xA}(T^{x})\eeqn
It will be shown that $F$ is an isomorphism, i.e.
\bq
F_{M}:=\oplus_{r \in \mtc{R}}(T_{2}^{xrx\inv, x})^{-1}_{M}(T_{2}^{x,r})_{M}:T^{x}(\ind_{B}^{A}(M))\ra \ind_{\;^xB}^{\;^xA}(T^{x}(M))
\eq
is an isomorphism in $\cc^{\;^{x}A}$ for any $M \in \cc^{B}$.  Indeed, one has to verify that the above morphism $F_{M}$ is compatible with the two equivariant structures defined above. This means that the following diagram:
{\Tiny
\begin{center}
\begin{tikzpicture}
  \matrix (m) [matrix of math nodes, row sep=6em,column sep=4em,minimum width=5em]
  { T^{xax^{-1}}(T^{x}\ind^{A}_{B}(M)) & T^{xax^{-1}}( \ind^{\;^{x}A}_{\;^{x}B}T^{x}(M)) \\
    T^{x}(\ind^{A}_{B}(M)) & \ind^{\;^{x}A}_{\;^{x}B}T^{x}(M)\\
    };%5
  \path[-stealth]
  (m-1-1) edge node [above] {$T^{xax^{-1}}(F_{M})$} (m-1-2)
    (m-2-1) edge node [above] {$F_{M}$} (m-2-2)
         
               (m-1-1) edge node [left] {$(\;^{x}\nu^{xax^{-1}})_{M}$} (m-2-1)
    (m-1-2) edge node [right] {$\eta^{xax^{-1}}_{M}$} (m-2-2)
   ;
\end{tikzpicture}
   \end{center}
}
is commutative.  On the components the above diagram becomes the following: \begin{center}
{\tiny
\begin{tikzpicture}
  \matrix (m) [matrix of math nodes, row sep=3.25 em,column sep=0.00005 em,minimum width=0.1em]
  { 
     T^{xax\inv}(T^{x}(T^{r}(M))) &  & T^{xax\inv}(T^{xr}(M)) &  & T^{xax\inv}(T^{xrx^{-1}}(T^{x}(M))) \\%1
     & (1) & &\;\; \;\;(2) & \\%2
    T^{xa} (T^{r}(M)) & & T^{xar}(M) =T^{xr'b}(M)  &  & T^{xarx\inv}(T^{x}(M))=T^{xr'bx\inv}(T^{x}(M))\\%3
     & (3) & & (4) & \\%4
     T^{x}(T^{a}(T^{r}(M)))& &  &  & T^{xr'x^{-1}}(T^{xbx^{-1}}(T^{x}(M)))\\%5 
    % & (4) & & (6) &  \\%6
    T^{x}(T^{ar}(M))=T^{x}(T^{r'b}(M)) & (5) & & (6)   & T^{xr'x^{-1}}(T^{xb}(M)) \\%6
   % & () & & (2) & \\%7
     T^{x}(T^{r'}(T^{b}(M)))) & & T^{xr'}(T^{b}(M)) & &  T^{xr'x^{-1}}(T^{x}(T^{b}(M))) \\%8
      T^{x}(T^{r'}(M))) & & T^{xr'}(M) & &  T^{xr'x^{-1}}(T^{x}(M)) \\%9
     };
  \path[-stealth]
    (m-1-1) edge node [left] {$(T_{2}^{xax\inv, x})_{T^{r}(M)}$} (m-3-1)%vertical A8
    %(m-1-3) edge [dashed] node [right] {$(T_{2}^{x, a})^{-1}_{T^{r}(M)}$} (m-5-1)%verticalA10
    (m-1-1) edge node [above] {$T^{xax\inv}(T_{2}^{x , r})_{M}$}(m-1-3)%A1horiz
  (m-3-5) edge node [right] {$(T_{2}^{xr'x\inv, xbx^{-1}})^{-1}_{T^{x}(M)}$} (m-5-5)%verticalA5
  (m-1-3) edge node [right] [above]{$(T_{2}^{x, a})^{-1}_{T^{r}(M)}$} (m-1-5)%horizA2
         (m-1-5) edge node [right] {$(T_{2}^{xax\inv, xrx^{-1}})_{T^{x}(M)}$} (m-3-5)%horizA3
          %(m-1-5) edge [dashed] node [left] {$(T_{2}^{xhx\inv, x})_{T^{l}(V)}$}
          (m-5-3)%horizA4
          % (m-5-3) edge [dashed] node [below] {$T^{xh}(\mu^{V}_{l})$}(m-5-5)%horiz A6
            (m-7-3) edge [dashed] node [right] {$(T_{2}^{xr' b})_{M}$}(m-3-3)%A14
            (m-7-1) edge [dashed] node [below] {$(T_{2}^{x,r'})_{T^{b}(M)})$} (m-7-3)%A13
             (m-1-3) edge [dashed] node [right] {$(T_{2}^{xax^{-1},xr})_{M}$} (m-3-3)%A13
         % (m-5-1) edge [dashed] node [below] {($T_{2}^{xh, l})^{-1}_{V}$} (m-5-3)%horizA11
    (m-3-1) edge node [left] {$(T_{2}^{x, a})^{-1}_{T^{r}(M)}$} (m-5-1)%verticalA9
    %(m-5-5) edge node [right] {$(T_{2}^{x,h})_{V}$} (m-9-5)%verticalA7
     %verticaql
      (m-7-3) edge [dashed]  node [right] {$T^{xr'}(\mu^{b}_{M})$} (m-8-3)
      (m-7-5) edge node [right] {$T^{xr'x^{-1}}(T^{x}(\mu^{b}_{M}))$} (m-8-5)
(m-6-5) edge node [left] {$T^{xr'x^{-1}}((T_{2}^{x, b})^{-1}_{M})$} (m-7-5)   
(m-5-5) edge node [right] {$T^{xr'x^{-1}}((T_{2}^{xbx^{-1}, x})_{M})$} (m-6-5)
(m-5-1) edge node [left] {$T^{x}((T_{2}^{a, r})^{-1}_{M})$} (m-6-1)%A12
(m-6-1) edge node [left] {$T^{x}((T_{2}^{r', b})^{-1}_{M})$} (m-7-1)%A12
    %(m-3-1) edge [dashed] node [right] {$\simeq$} (m-2-2)
  (m-7-3) edge [dashed] node [below] {$(T^{xr'x^{-1}, x}_{2})^{-1}_{T^{b}(M)}$} (m-7-5)%A15
    (m-8-1) edge node [below] {$(T^{x, r'}_{2})_{T^{b}(M)}$} (m-8-3)%A15
(m-8-3) edge node [below] {$(T^{xr'x^{-1}, x}_{2})^{-1}_{M}$} (m-8-5)%A15
(m-3-1) edge [dashed] node [below] {$(T^{xa, r}_{2})_{M}$} (m-3-3)%A15
(m-3-3) edge [dashed] node [below] {$(T^{xr'x^{-1}, x}_{2})^{-1}_{T^{b}(M)}$} (m-3-5)%A15
(m-3-3) edge [dashed] node [left] [above] {$(T^{xr'x^{-1}, xb}_{2})^{-1}_{M}$} (m-6-5)%A15
(m-3-3) edge [dashed] node [below] {$(T^{x, r'b}_{2})^{-1}_{T^{b}(M)}$} (m-6-1)%A15
   (m-7-1) edge node [left] {$T^{x}(T^{r'}(\mu^{b}_{M}))$} (m-8-1)%A16
      %(m-9-5) edge node [below] {$T^{x}(\mu^{h}_{V})$} (m-9-1)%A17
         ;
\end{tikzpicture}
}
\end{center}
Note that diagrams $(1)-(6)$ are commutativity by the associativity of the action, Equation \eqref{ro-2}. The bottom two rectangles are commutative since $T^{x, r'}_{2}$ and $T^{xr'x^{-1}, x}_{2}$ are natural transformations.

\epf
\ncm{\wta}{With the above notations one has}
\bl\label{tensorconj}
Let $G$ be a finite group acting by tensor autoequivalences on the tensor category $\cc$ and $A$ be a subgroup of $G$. \wta  
\bq
T^{x}(M\ot N)\simeq T^{x}(M)\ot T^{x}(N)
\eq
for any two objects  $M, N \in \cc^{\;^{x}A}$. Therefore $T^{x}$ is a tensor isomorphism between $\cc^{A}$ and $\cc^{\;^{x}A}$.
\el
\bpf 
One has to check that the tensor structure $(T^{x}_{2})^{M, N} :T^{x}(M\ot N)\ra T^{x}(M)\ot T^{x}(N)$ of $T^{x}$ is a morphism in $\cc^{\;^{xA}}$. Thus for any $M, N \in \cc^{\;^{x}A}$ one has to check the commutativity of the following digram:

{\SMALL
\begin{tikzpicture}
  \matrix (m) [matrix of math nodes, row sep=2em,column sep= 0.9 em,minimum width=3em]
  { T^{xax\inv}(T^{x}(M\ot N)) &   &  T^{xax\inv}(T^{x}(M)\ot T^{x}(N))\\%1
  & (1) & \\ %2\\
   &  & T^{xax\inv}(T^{x}(M))\ot T^{xax\inv}(T^{x}(N))\\%3
   T^{xa}(M\ot N) & & T^{xa}(M)\ot T^{xa}(N)\\%4
    & (2) & \\%5
  T^{x}(T^{a}(M\ot N)) & & \\%6
  %T^{xa}(M)\ot T^{xa}(N)\\%7
  T^{x}(T^{a}(M)\ot T^{a}(N)) & & T^{x}(T^{a}(M))\ot T^{x}(T^{a}(N))  \\%7
  & (3) & \\%8\\
  T^{x}(M\ot N) & &T^{x}(M)\ot T^{x}(N)\\};%9
  \path[-stealth]
  (m-1-1) edge node [above] {$T^{xax\inv}((T^{x}_{2})^{M, N})$} (m-1-3)
    (m-9-1) edge node [above] {${(T^{x}_{2})}^{M, N}$} (m-9-3)
         
         (m-4-1) edge node [left] {$(T^{x,a}_{2})^{-1}_{M\ot N}$} (m-6-1)
         (m-6-1) edge node [left] {$T^{x}((T_{2}^{a})^{M, N})$} (m-7-1)
         (m-7-1) edge node [left] {$T^{x}(\mu^{M}_{a}\ot \mu^{N}_{a})$} (m-9-1)
         
            (m-3-3) edge node [right] {$(T_{2}^{xax\inv, x})_{M}\ot (T_{2}^{xax\inv, x})_{N}$} (m-4-3)
         (m-4-3) edge node [right] {$(T_{2}^{x, a})^{{-1}}_{M}\ot (T_{2}^{x, a})^{{-1}}_{N}$} (m-7-3)
         (m-7-3) edge node [right] {$ T^{x}(\mu^{M}_{a}) \ot T^{x}(\mu^{N}_{a})$} (m-9-3)
               (m-1-1) edge node [left] {$(T^{xax\inv, x}_{2})_{M\ot N}$} (m-4-1)
    (m-1-3) edge node [right] {$(T^{xax\inv}_{2})^{{T^{x}(M), T^{x}(N)}}$} (m-3-3)
    (m-4-1) edge [dashed] node [above] {$(T^{xa}_{2})^{M,N}$} (m-4-3)
     (m-7-1) edge [dashed] node [above] {$({T_{2}^{x})}^{{T^{a}(M)},{T^{a}(N)}}$} (m-7-3)
    ;
\end{tikzpicture}
}

The upper pentagon (1) is commutative since $T_{2}^{xax\inv , \;x}$ is a natural transformation of tensor functors, Equation \eqref{tensor-rho}. The middle pentagon $(2)$ is commutative since $T^{x, a}_{2}$ is a natural transformation of tensor functors, same Equation \eqref{tensor-rho}. The bottom rectangle commutes from the compatibility condition of the tensor functor $T^{x}$ with the tensor product of morphisms, Equation \eqref{natf}.\epf

\subsection{$ G$-functors} Let $G$ be a finite group.
A Mackey functor (or a $G$-functor) over a ring $R$ is a collection of $R$-modules $\{a(H)\}_{H\leq G}$
together 
with morphisms
$I^{H}_{K} : a(K)\ra a(H)$,
$R^{H}_{K}: a(H)\ra a(K)$ and
$c_{H, g} : a(H)\ra a( \;^{g}H)$
for all subgroups $H$ and $K$ of $G$ with $K\leq H$ and for all $g \in G$. This datum satisfies the following compatibility conditions:

$\mathrm{(M0)}$ $I_{H}^{H}, \;R_{H}^{H}, c_{H, h} : M(H)\ra M(H)$
 are the identity morphisms for all subgroups H and $h \in H$.

$\mathrm{(M1)}$ $R_{K}^{J}R_{H}^{K}=R^{J}_{K}$ for all subgroups $J\leq K\leq H$.

$\mathrm{(M2)}$ $
I_{H}^{K}I_{J}^{H}=I^{K}_{J}$, for all subgroups $J\leq K\leq H$.

$\mathrm{(M3)}$
$c_{H, g}c_{\;^{g}H, h} = c_{H, gh} $, for all $H \leq G$ and $g, h \in G$.

$\mathrm{(M4)}$
For any subgroups $K, H \leq G$ the following Mackey relation is satisfied:
\beq\label{mackeyss}
R^{G}_{H}I^{G}_{K}=\bigoplus_{x \in H\backslash G\slash K}I^{H}_{\;^{x}K\cap H}R^{^{x}K}_{^{x}K\cap H}c_{K, x}.
\eeq
Moreover, a Green functor over a commutative ring $R$, is a $G$-functor $a$
such that for any subgroup $H$ of $G$ one has that 
$a(H)$ is an associative $R$-algebra with identity and satisfying the following:

$\mathrm{(G1)}$ $R_{H}^{K}\;\;
\text{and} \;\; c_{{H, g}}\;\;\text{are always unitary R-algebra homomorphisms,}
$

$\mathrm{(G2)}$ $
I_{H}^{K}(aR^{H}_{K}(b)) = I_{H}^{K}(a)b
$,

$\mathrm{(G3)}$ $
I_{H}^{K}(R^{H}_{K}(b)a) =b I_{H}^{K}(a)
$
 for all subgroups $K$ and $H$ and all $a \in a(K)$ and $b \in b(H)$.
 %\subsection{$K$-theory for exact categories}
 %In \cite{quill} Quillen associated to an abelian (exact) category $\ca$ a topological space
%$BQA$ such that exact functors between categories define the continuous maps between the corresponding spaces and isomorphisms of functors define the homotopies between the corresponding maps. In other words, QuillenÕs space is a
%2-functor from the 2-category of abelian (exact) categories (with isomorphisms
%of functors as 2-morphisms) to the 2-category of topological spaces.
%Waldhausen [27] proved that the 2-functor K = 
%BQ is permutable (in some
%sense) with the product. Namely, he constructed the continuous map K(A) ?
%K(B) ? K(C) for any biexact functor A ? B ? C.
%The homotopy groups $K_{*}(A) = \pi_{*}(A)$ of the Waldhausen space K(A) are called algebraic $K$-theory of the category $A$. 
\subsection{Proof of Theorem \ref{main2}}
\bpf
Similarly to \cite{kevin} we use the following elementary facts about K-theory (see \cite{quillen}): If
$F_{1}$ and $F_{2}$ are isomorphic exact functors on an exact category, then they induce the
same map on $K$-theory;  and if $F_{1}$ and $F_{2}$ are exact functors on an exact category inducing
homomorphisms $f_{1}$ and $f_{2}$ on $K$-groups, then the functor $F_{1}\oplus F_{2}$ induces the
homomorphism $f_{1}+f_{2}$.
Now identities $\mathrm{(M1)}$ - $\mathrm{(M4)}$ follow from their functorial counterpart proven in the previous section.

%Clearly the induced maps $c_{H, h}:K_{i}(\cc^{H})\ra K_{i}(\cc^{H})$ are identity maps since one has $T^{h}(M)\simeq M$ for any $M\in \cc^{H}$. 
For example, the identity $\mathrm{(M2)}$ follows from the equality $\ind_{K}^{H}\ind_{J}^{K}=\ind_{J}^{H}$ which can be verified by a straightforward computation.
\epf
\subsection{Proof of Theorem \ref{grf}} First we need the following proposition:
\bp\lb{atp}Suppose that $\cc, \cd$ and $ \ce$ are rigid monoidal categories and $F_1:\cc \ra \cd$ and $F_2:\cc \ra \ce$ are two monoidal functors with left adjoint functors $I_1:\cd\ra \cc$ and respectively $I_2:\ce\ra \cc$.  Then for any objects $M \in \co(\cd)$ and $N \in \co(\ce)$ one has the canonical isomorphism in $\cc$
\beq\lb{ti}
I_1(M)\ot I_2(N)\simeq I_1(F_{1}(I_2(N))\ot M) \simeq I_2(F_2(I_1(M))\ot N).
\eeq

\ep
\bpf
It can be shown by a straightforward computation that 
\beq\lb{ti}
\Hom_{\cc}(I_1(M)\ot I_2(N), P)\simeq \Hom_{\cc}(I_2(F_2(I_1(M))\ot N)
, P)
\eeq
for any object $P \in \cc$. Indeed, 
\begin{eqnarray*}
  \Hom_{\C}(I_1(M)\ot I_2(N), P)&=& \Hom_{\C}(I_{1}(M),\;P \ot I_{2}(N)^*)= \\
  =\Hom_{\mtc{D}}(M,\;F_{1}(I_{2}(N)^*\ot P)) &=&  \Hom_{\mtc{D}}(M, F_{1}(I_{2}(N))^* \ot F_{1}(P))=\\
  =\Hom_{\mtc{D}}(F_{1}(I_{2}(N)) \ot M, F(P)) &=&  \Hom_{\mtc{C}}(I_1(F_{1}(I_2(N))\ot M), P)
\end{eqnarray*}

Then Yoneda's lemma  implies the conclusion.\epf
%\blue{closed categories bvl}
%\subsection{On the Grothendieck ring of an equivaraintized fusion category}
%Put $A(g):=K_{0}(\cc^{F_{g}})$ and $A=K_{0}(\cc^{F})$

In particular for $\ce=\cc$ and $F_{2}=I_{2}=\id_{\cc}$  one obtains that 
\beq\lb{ti2}
I_{1}(M)\ot V\simeq I_{1}(M\ot F_{1}(V))
\eeq for any objects $M\in \cd$ and $V\in \cc$.
\subsubsection{Proof of Theorem \ref{grf}}

Note that the multiplicative relations $\mathrm{(G2)}$ and $\mathrm{(G3)}$ follow from Proposition \ref{atp}.

\bn{example}
Suppose that $R \subset S$ is a Galois extension of rings with Galois Group $G$. Then as before $G$ acts on the category $S$-mod and $(S$-mod $)^{G}\simeq S\# \Z G$-mod. Since $S^{G}$ is Morita equivalent to $S\# \Z G$ and the $K$ -theory is preserved by Morita equivalence it follows that our results extends the results from \cite{kevin}.% \blue{ work with abelian case, to fully recover the result}
\end{example}
%\bibliographystyle{amsplain}
%\bibliography{bob}
%\ed
%\np
\section{Coherent group actions on graded  fusion categories}\lb{hmg} Let $\cc$ be a graded fusion category by a finite group $G$. Recall that this means that
$
\cc=\oplus_{g \in G}\cc_{g}
$
as abelian categories, and the tensor functor $\ot: \cc \times \cc\ra \cc $ sends $\cc_g\ot \cc_h$ into $\cc_{gh}$. For an object $V \in \cc$ define by $V_g$ the homogenous component of $V$ of degree $g$ from the above grading.% of $\cc$.
 \md
 
Suppose further that another finite group $F$ acts by group automorphisms on $G$. Suppose that $F$ also acts by tensor automorphisms on the category $\cc$ via the action $T:F\ra \underline{\mtr{Aut}}_{\ot}(\cc)$ given by $x \mto T^x  :\cc \ra \cc$.

\bn{defn} An action of a finite group $F$ on the $G$-graded fusion category $\cc$ is called {\it coherent} with respect to the action of $F$ on $G$
 if 
\beq
T^x(\cc_g)\subset \cc_{\;^xg}
\eeq
for all $x \in F$ and $g \in G$.
\end{defn}
\md
\bp Suppose that a finite group $F$ acts by tensor automorphisms on a $G$-graded fusion category $\cc$ via the action $T:F\ra \underline{\mtr{Aut}}_{\ot}(\cc)$ given by $x \mto T^x  :\cc \ra \cc$. Suppose further that $\cc_1$ is stable under the action of $F$. Then there is an action of $F$ on $G$ by group automorphisms such that the action of $F$ on $\cc$ is coherent with respect to this action.  
\ep
\bpf 
Let $V, W\in \cc_g$ and suppose that $T^x(V)\in \cc_{g_1}$ and $T^x(W) \in \cc_{g_2}$. Then $T^x(V^*\ot W)\in \cc_1$ since $V^*\ot W \in \cc_1$. On the other hand $T^x(V^*\ot W)\simeq T^x(V^*)\ot T^x(W)\in \cc_{{g_1}^{-1}g_2}$ which implies that $g_1=g_2$. Denoting $\;^xg:=g_1$ then it is easy to check that this defines an action of $F$ on $G$ by group automorphisms.
\epf
\subsection{Examples of coherent actions of groups and their equivariantized categories} In this subsection we give some examples of coherent group actions on fusion categories.
\bn{example}\lb{gbb}{\it Braided $G$-crossed categories.}
Recall \cite{tuv} that a braided $G$-crossed fusion category is a quadruple $(\cc, G, T, c)$, where $G$ is a finite group, $\cc$ is a tensor category with a (not necessarily faithful) $G$-grading $\cc=\oplus_{{g \in G}}\cc_{g}$
 and a tensor action $T : G \ra \underline{\mtr{Aut}}_{\ot}(\cc)$, $g \mapsto T_{g}$ satisfying $T_{g}(\cc_{h})\subseteq \cc_{ghg\inv}$. Moreover the crossed braiding $c$ is defined by $c(X, Y):X\ot Y \ra T_{g}(Y)\ot X$ for all $X\in \cc_{g}$ and $Y \in \cc$. The compatibility conditions that have to be satisfied by  this datum can be found for example in \cite{tuv, DGNO}. Note that in this case $F=G$ acts coherently on $\cc$ with respect to the action of $G$ on  itself given by conjugation.\end{example}
 \bn{example}\lb{twisted} 
Group actions on graded pointed fusion categories are always coherent. See \cite[Section 4]{naidu}. In particular,  it follows by \cite[Lemma 6.3]{naidu} that the representation category of a (twisted) quantum double of a finite group is the equivariantization of a coherent action.
\end{example}
 \subsubsection{Universal gradings for the category of representations of semisimple Hopf algebras}
%\blue{Recall the grading of $\rep(B)$ and eliminate all the other identities.}
Let $A$ be a semisimple Hopf algebra over an algebraically closed field $k$. It is well known that $\mtr{Rep}(A)$ is a fusion category. Moreover there is a maximal central Hopf subalgebra $K(A)$ of $A$ such that $\mtr{Rep}(A//K(A))$ coincides to $\mtr{Rep}(A)_{
_{ad}}$ the adjoint subcategory of $\rep(A)$, see \cite[Theorem 2.4]{GN}. Since $K(A)$ is commutative it follows that
$K(A)=kG^*$ where $G$ is the universal grading group of $\mtr{Rep}(A)$. The (universal) grading on $\rep(A)$ is given by  \bq\label{grdm}
\rep(A)_{g}=\{M \in \Irr(A)\;|\; p_{g}m=m \;\text{for all}\;m\in M\}.
\eq
Here $p_{g}\in k^{G}$ is the dual basis of group element basis $g \in G$.

\subsubsection{Cocentral extensions of semisimple Hopf algebras}\label{cocs}

%\red{One can also put it as an appendix}
Suppose that we have a cocentral extension of Hopf algebras 
\beq\lb{cocen}
k \ra B\xrightarrow{i} H\xrightarrow{\pi} kF\ra k.
\eeq

Recall that this means the above sequence is exact \cite{schss} and  $kF^*\subset \mtc{Z}(H^*)$ via $\pi^*$. 
On the other hand, using the reconstruction theorem from \cite{AD} it follows that 
$
H \simeq B\;^{\tau}\#_{\sg} \;kF
$
for some cocycle $\sg:B\ot B \ra kF$ and some dual cocycle $\tau:kF \ra B\ot B$ satisfying certain compatibility axioms.  In this case there is a weak action of $F$ on $B$ denoted by $f.b$ such that the multiplication and comultiplication on $H$ become
\bq\lb{mu}
(b\#_{\sg}f)(c\#_{\sg}g)=b(f.c)\sg(f,g)\#_{\sg}fg
\eq
and respectively
%\beq\Delta(b\#g)=(b_1\tau(g)_j\#_{\sg}g)\ot (b_2\tau(g)^j\#_{\sg}g)\eeq
\bq\lb{com}
\D(b\#_{\sg}\bar{f})=(b_1\tau(\bar{f})_i\#_{\sg}\bar{f})\ot(b_2\tau(\bar{f})^i\#_{\sg}\bar{f}).
\eq
%Moreover the antipode is given by the formula $$S(a\#_{\sg}g)=g^{-1}S(a)=(g^{-1}.S(a))\#_{\sg}g^{-1}.$$
%and the weak cocation $T$ is denoted by $T(b)=\sum_i b_i\ot b^i \in B \ot kF$.
\bp\label{cocsp} Suppose that we have a cocentral extension of semisimple Hopf algebras as in Equation \eqref{cocen}.
%\begin{equation*}
%k \ra B\xrightarrow{i} H\xrightarrow{\pi} kF\ra k.%\end{equation*}
Moreover, with the above notations suppose that $\cc:=\rep(B)$ and let
$\mtc{C}=\oplus_{g \in G}\mtc{C}_g
$
be the universal grading of $\rep(B)$ where $K(B)=kG^*$. Then $F$ acts coherently on $\cc$ with respect to a given action by group automorphisms of $F$ on $G$.
\ep
\bpf
First it will be shown that $F$ acts on $G$ by group automorphisms. In order to do this we show that $F$ acts on $K(B)=kG^*$ by Hopf automorphisms.
\md Note that if $b \in K(B)$ then $f.b \in K(B)$ where $``.``$ represents the weak action of $F$ on $B$ from above. Indeed $f.b \in \cz(B)$ since $F$ acts by algebra automorhisms on $B$.
On the other hand using formula $(A)$ from \cite[Section 2]{AD} it follows that 
\beq
\Delta(f.b)=\tau(f)(f.b_1\ot f.b_2)\tau(f)^{-1}.
\eeq
Therefore if $b \in \cz(B)$ then 
$
\Delta(f.b)=f.b_1\ot f.b_2
$
which shows that $F. K(B)$ is a central Hopf subalgebra of $B$.  Thus $F.K(B)\subseteq K(B)$ and $F$ acts by Hopf algebra automorphisms on $K(B)$. This implies that there is an action of $F$ on $G$ such that the action of $F$ on $K(B)=kG^*$ is given by
$
x.p_g=p_{\;^xg}
$
for all $x \in F$ and $g \in G$.
\md Following \cite[Proposition 3.5]{natalecoc} it follows that $F$ acts on the fusion category $\Rep(B)$ and $\Rep(H)=\Rep(B)^F$.  Recall from \cite{natalecoc} that for all $x \in F$ the action is given by $T^x(M)=M$ as vector spaces with the action of $B$ on $T^x(M)$ given by
$
b.\;^xm=(x^{-1}.b)m.
$

\md In order to verify  that the above action of $F$ on $\rep(B)$ is coherent with respect to  this action one has to verify that if $M \in \co(\cc_g)$ then $T^x(M)\in \co( \cc_{\;^xg})$. Using Equation \eqref{grdm} it follows that for any $h \in G$ one has that 
%\begin{equation*}
${p_h.}\;^xm=(x^{-1}.p_h)m=p_{\;^{x^{-1}}h}\;m=\delta_{\;^{x^{-1}}h,\; g}\;m=\delta_{h, \;^{x}g}\;m$
%\end{equation*}
which shows that indeed $T^x(M)\in \co(\cc_{\;^xg})$.
\epf

%\br If $k \ra k^{G}\xrightarrow{i} H\xrightarrow{\pi} kF\ra k $
%is an abelian cocentral extension then the action of $F$ on $K(k^{G})=k^{G}$ from above Theorem corresponds with the usual action of $F$ on $k^{G}$ from \cite{Masext}.
%\er

\sub{Simple objects for equivariantizations of coherent actions}In this subsection we investigate the simple objects of an equivariantization under a coherent action. 
\md

Suppose that $F$ acts coherently on a $G$-graded fusion category $\cc$ with respect to a given action of $F$ on $G$.
With the above notations note that the stabilizer $F_{g}$ of an element $g \in G$ acts by $k$-linear automorphisms on the abelian subcategory $\cc_{g}$ of $\cc$. In particular one obtains in this way an action of  $F$ on the fusion subcategory $\cc_{1}$ of $\cc$. Note that this action on $\cc_{1}$ is by tensor automorphisms.

\bl With the above notations, suppose that $(V, {\mu_V^x}|_{\{x \in F\}})\in \cc^F$ is an equivariantized object with a canonical decomposition $V=\oplus_{g \in G}V_{g}$. Then for all $g \in G$ one has that $(V_g,  {\mu^x_{V_{g}}}|_{\{x \in F_g\}}) \in \cc^{F_g}$.
\el

\bpf If $V=\oplus_{g \in G}V_{g}$ then $T^{x}(V)=\oplus_{g \in G}T^{x}(V_{g})$.
Since $\mu_V^x:T^x(V)\ra V$ is an isomorphism in $\cc$ it follows that $\mu^x_V$ sends the component $T^{x}(V_g)$ of $T^{x}(V)$ to the component $V_{\;^{x}g}$ of $V$. Thus if $x \in F_g$ then $\mu^{x}_{V}|_{T^{x}(V_{g})}:T^{x}(V_g)\ra V_{g}$ is an isomorphism. Note that the compatibility conditions from Equation \eqref{deltau} for an equivariantized object $V \in \cc^{G}$ implies that  $V_{g} \in \cc^{F_{g}}$. \epf

Let $\Gm$ be a set of representatives for the orbits of the action of $F$ on the group $G$. For any $g \in G$ let $\co(g)$ be the orbit of $g$ under the action of $F$. Next proposition is a generalization of \cite[Proposition 2.7]{gnn}.

\bp\label{desgnn} Suppose that $F$ acts coherently on a  $G$-graded fusion category $\cc$ with respect to an action by group automorphisms on $G$. Then the set of simple objects of $\cc^F$ is parametrized by pairs $(g , M)$ where $g \in \Gm$ and $M \in {\cc_g}^{F_g}$  is a simple object. The simple object associated to the pair $(g, M)$ is given by the induced object $\Ind_{F_{g}}^{F}(M)$.
\ep 

\bpf  
Suppose that $(V, {\mu^{x}_{V}}_{\{x \in F\}})\in \cc^F$ is a simple object and let $V=\oplus_{g \in G}V_g$ with $V_{g }\in \cc_{g}$ be its decomposition viewed as an object of $\cc$. If $V_g \neq 0$ then $T^x(V_g)=V_{\;^{x}g}$ is also not zero and therefore $\opl_{h \in \co(g)}V_h$ is an equivariantized object of $\cc^F$. Since $V$ is a simple object it follows that $V$ is supported only on one orbit, namely the orbit $\co(g)$ of $g$. Thus $V=\oplus_{{h \in \co(g)}}V_{h}$. Moreover it follows that $V_g\in \cc_g^{F_g}$ and one can associate to the object $(V, {\mu^{x}_{V}}_{\{x \in F\}})\in \cc^F$ the pair $(g, M)$ with $M:=V_g$. Clearly $M$ is a simple object of $\cc_{g}^{F_{g}}$ if $V$ is a simple  object of $\cc^{F}$.
\md
Conversely, to any pair $(g, M)$ as above one associates the simple object $\Ind_{F_{g}}^{F}(M) \in \cc^{G}$.
\md
We have to show that the above two constructions are one inverse to the other.  It is easy to see that the component of order $g$ of $\Ind_{F_{g}}^{F}(M) \in \cc^{G}$ coincides to $M$ as an object of $\cc^{F_{g}}_{g}$.
\md Thus it remains to show that $V\simeq \ind^{F}_{F_{g}}(V_{g})$ if $V=\oplus_{g \in G}V_{g}$ is a simple object of $\cc^{F}$. 
This is equivalent with showing that
\bq
V\simeq \oplus_{r \in F/F_{g}}T^{r}(V_{g})
\eq
as objects of $\cc^{F}$.

Since $\mu_{V}^{x}:T^{x}(V) \xra{\simeq} V$ it follows that $\mu_{V}^{x}|_{T^{x}(V_{g})}:T^{x}(V_{g}) \xra{\simeq}V_{\;^{xg}}$ are also isomorphisms in $\cc$. This shows that $f: \ind_{F_{g}}^{F}(V_{g})\xra{\simeq} V$ is an isomorphism in $\cc^{F}$ where $f:=\oplus_{r \in F/F_{g}}\mu_{V}^{r}|_{T^{r}(V_{g})}$.
\epf
\br
Note that there is an embedding $\cc_{g}^{F_{g}}\subset \cc^{F_{g}}$ and one can apply the induced factor $\ind^{F}_{F_{g}}$ to $M$.
\er
%}
%\bc Suppose that $F$ acts homogeneously on a graded fusion category $\cc$ and
%let $V=\oplus_{h \in \co(g)}V_{h}$ be a simple
%object of $\cc^{F}$. Then 
%\bq\Ind^{F}_{F_{g}}(V_{g})\simeq \ind_{F_{h}}^{F}(V_{h})\eq
%for all $h \in \co(g)$.\ec
%\bpfBoth terms are isomorphic to $V$.\epf\md
\subsection{Definition of $S(g, M)$} Let $M \in \cc^{F_{g}}$ be a simple object.
Define by $S(g, M)$ as the simple induced object $\ind^{F}_{F_{g}}(M) \in \cc^F$ from above. Thus as objects of $\cc$ one has that $S(g, M):=\oplus_{r \in F/F_g}T^r(M)$ with the equivariant $F$-structure obtained from Equation \eqref{muind}.

\br\lb{shift} Note that if $(M , \mu_{M}) \in \cc_g^{F_g}$ then $(T^x(M), \;^x\mu_{M})\in \cc_{\;^xg}^{F_{\;^xg}}$ by Proposition \ref{conj}. The proof of previous theorem implies that:
\beq
S(g, M)\simeq S(\;^xg, T^x(M))
\eeq
for any $x \in F$, $g \in G$ and $M \in \cc_g^{F_g}$.
%\subsection{Mackey's decomposition of an induced object}
\er
\subsection{Tensor product formula}
Let $M\in \cc_g^{F_{g}}$ and $N \in \cc_h^{F_{h}}$ be two equivariant objects. Define the following equivariant object
\beq\label{mgh}
m_{g,h}(M, N):=\ind_{F_{g}\cap F_{h}}^{F_{gh}}(\res^{F_{g}}_{F_{g}\cap F_h}(M)\ot \res_{F_{g}\cap F_h}^{F_h}(N))\in \cc_{gh}^{F_{gh}}.
\eeq
%\ind_{F_{\;^xg}\cap F_{h}}^{F_{\;^xgh}}(\res^{F_{\;^xg}}_{F_{\;^xg}\cap F_h}(\;^xM)\ot \res_{F_{\;^xg}\cap F_h}^{F_h}(N))
\bt \lb{tpgr} Suppose that $F$ acts coherently on the $G$-graded fusion category $\cc$ with respect to a given action by group automorphisms of $F$ on $G$. With the above notations one has
\beq\lb{tepg}
S({g, M})\ot S({h,N})\simeq \opl_{x \in D}S(\;^xgh, \;\;m_{ _{\;^xg, h}}(T^x(M),N))
\eeq
where $D$ is a set of representatives for the double cosets $F_h \backslash F / F_g$.
\et
\bpf
One has by definition $S(g, M)=\ind_{F_g}^F(M)$. Applying formula (\ref{ti}) one has that
\beqarn
S({g, M})\ot S({h,N}) & = & \ind_{F_g}^F(M)\ot \ind_{F_h}^F(N)\\
& \simeq & \ind_{F_h}^F(\res_{F_h}^F(\ind_{F_g}^F(M)\ot N)).
\eeqarn

On the other hand applying Theorem \ref{macky} one has that 
\beqarn
\res_{F_h}^F(\ind_{F_g}^F(M))\simeq \oplus_{x \in D}\ind_{{F_h} \cap \;^x{F_g}}^{F_h}(\res^{\;^x{F_g}}_{\;^x{F_g}\cap {F_h}}(T^x(M)))
\eeqarn
Then applying Equation \eqref{ti2} one obtains that
\beqarn
\res_{F_h}^F(\ind_{F_g}^F(M))\ot N\simeq \oplus_{x \in D}(\ind_{{F_h} \cap \;^x{F_g}}^{F_h}(\res^{\;^x{F_g}}_{\;^x{F_g}\cap {F_h}}(T^x(M))\ot N) \\ \simeq \oplus_{x \in D}(\ind_{{F_h} \cap \;^x{F_g}}^{F_h}(\res^{\;^x{F_g}}_{\;^x{F_g}\cap {F_h}}(T^x(M)\ot \res^{{F_h}}_{\;^x{F_g}\cap {F_h}}(N))
\eeqarn
and therefore 
\beqn S({g, M})\ot S({h,N}) 
\simeq  \oplus_{x \in D}\ind_{\;^{x}F_{g}\cap F_{h}}^{F}(\res^{\;^{x}F_{g}}_{\;^{x}F_{g}\cap F_h}(T^{x}(M))\ot \res_{\;^{x}F_{g}\cap F_h}^{F_h}(N))
\eeqn
which by definition coincides to $\opl_{x \in D}S(\;^xgh, \;\;m_{ _{\;^xg, h}}(T^x(M),N))$.
\epf
\bc Suppose that $F$ acts coherently on the $G$-graded fusion category $\cc$ with respect to a given action by group automorphisms of $F$ on $G$. With the above notations it follows that
\beq\lb{ditf}
\cc^F \simeq \oplus_{g \in \Gm}\cc_g^{F_g}
\eeq
as indecomposable $\cc_{1}^F$-bimodule categories, where $\Gm$ is a set of representative elements for the orbits of the action of $F$ on $G$ .\ec
\bpf
Remark that $\cc_{1}^{F}$ is a tensor subcategory of $\cc^{F}$ consisting on those objects of $\cc^{F}$ supported only on $\cc_{1}$. Define the functor $F: \cc^{F} \ra \oplus_{g \in \Gm}\cc_g^{F_g}$ by sending $S(g, M)\mapsto M$. Note that formula (\ref{tepg}) implies that 
\beqn
S(1, M)\ot S(h, N)\simeq S(h, \res^{F}_{F_{h}}(M)\ot N)
\eeqn
and
\beqn
S(h, N)\ot S(1, M)\simeq S(h, N\ot \res^{F}_{F_{h}}(M))
\eeqn
which shows that each $\cc^{F_{h}}$ is a $\cc_{1}^{F}$-bimodule category. The above formulae also show that $\cc_{h}^{F_{h}}$ is an indecomposable $\cc_{1}^{F}$-bimodule category for all $h \in H$.
\epf
\subsection{On the Grothendieck ring of an equivariantization under a coherent action}\lb{withsp}In this subsection we show that the Grothendieck ring of an equivariantization under a coherent action has the structure of a Green ring as introduced in \cite{scoh}.
\subsubsection{On the rings introduced by Witherspoon and Bouc} In this subsection we recall the Green rings introduced in \cite{scoh}. 

Let $F$ be a finite group acting by group automoprphisms on another finite group $G$. Suppose  that 
$A=\oplus_{g \in G}A(g)$ is a graded vector space endowed with two linear structures: $m_{g,h}:A(g)\ot A(h)\ra A(gh)$
and  $c_{x,g}:A(g)\ra A(\;^xg)$ satisfying the following compatibilities:

$(\mathrm{C1})$
$c_1=\id \text{and\;} c_{xy}=c_x c_y$ where 
 $c_x:=\oplus_{g \in G}\;c_{x,g}:A\ra A$.

$(\mathrm{C2})$ $c_{x, g}=\id_{A(g)}$ if $x \in F_{g}$.

$(\mathrm{C3})$ $c_xm_{g,h}=m_{\;^xg,\;^xh}(c_x\times c_x)
$

$(\mathrm{C}4)$ There is an element $1 \in A(1)$ such that
$
c_x(1)=1 \; \text{for all\;} x \in L \; \text{and\;} m_{1,g}(1, \al_g)=m_{g,1}(\al_g,1)=\al_g
$.

Let $A^F:=\{a \in A\;|\; c_x(a)=a\;\text{for all}\;x\in F\}$
be the subspace of $F$-invariants elements of $A$.

$(\mathrm{C5})$ For any $g \in G$ and $\al, \beta, \gamma \in A^{F}$ one has that
\beqn
\sum_{(d,e,f)\in T_g}m_{de, f}(m_{d,e}(\al_d, \beta_e),\gm_f)=\sum_{(d,e,f)\in T_g}m_{d, ef}(\al_d, m_{e,f}(\beta_e,\gm_f))
\eeqn
where the set $T_g$ is defined as follows. Note that the stabilizer subgroup $F_g$ acts diagonally on the set $\{(d,e,f)\in G\times G\times G\;|\; def=g\}$. By $(\mathrm{C}1)$-$(\mathrm{C}4)$ the left and right members of the previous equality do not depend on the chosen set $T_g$ of representative elements for the orbits of the action of $F_g$ on the above set.

Define a multiplicative structure on $A^{F}$ given for $\al, \beta \in A^{F}$ by 
\beq\tag{M}\label{multo}
(\al\beta)_g=\sum_{\{(h,k)\in F_g/G \times G\;| \;hk=g\}}m_{h,k}(\al_h,\beta_k)
\eeq
where the action of $F_g$ on $G\times G$ is diagonal.
Then it is shown in \cite{scoh} that under the conditions $(\mathrm{C}1)$-$(\mathrm{C}4)$ the above multiplication on $A^F$ is associative if and only if condition $(\mathrm{C}5)$ is also satisfied. Moreover, this multiplication not depend on the choice of representative set $T_{g}$.% of  $F_g/G \times G$.%\subsubsection{}
%Fix a set of representatives $\{ g_1, g_2,\cdots, g_r \}$ for the orbits of the action of $F$ on $G$. Then one can identify
%\beq\lb{identific}
% A^F\simeq\oplus_{i=1}^rA(g_i)^{F_{g_i}} 
% \eeq
%and it was noticed in \cite{scoh} that under this identification the multiplication from Equation \eqref{multo}  becomes:\beq\lb{nm}\al_i\beta_j=\sum_{x \in D}
%m_{\;^{yx}g_i\;^yg_j}(c_{yx, g_i}(\al_i), c_{x, g_i}(\beta_j))\eeq
%where $y\in F$ is chosen such that $\;^{yx}g_i\;^yg_j=g_k$ for a uniquely determined $k=k(x)$.
%%%%%%%
\bt\lb{idr}
Let $F$ be a finite group acting coherently on a $G$-graded fusion category $\cc$ with respect to a given action by group automorphisms of $F$ on $G$. Then the Grothendieck ring of $\cc^G$ has the multiplicative structure of Equation \eqref{multo}.
\et
\bpf
Let $A=\oplus_{g \in G} A(g)$ where $A(g)=K_0(\cc^{F_g})$. Using Proposition \ref{conj} one can efine $c_{x, g}: A(g)\ra A(\;^xg)$ by $[M]\mapsto [T^x(M)]$. Define also $m_{g,h}:A(g)\times A(h) \ra A(gh)$ via the map $m_{g,h}$ from Equation \eqref{mgh}. Then it is easy to verify that the compatibility conditions $(\mathrm{C}1)$,  $(\mathrm{C2})$ and $(\mathrm{C4})$ from the previous subsection are satisfied. Condition $(\mathrm{C3})$ is verified in the lemma below.
Moreover it is clear that $K_0(\cc^F)\hookrightarrow A^F$ via $[S(g, M)] \mapsto \oplus_{r \in F/F_{g}}[T^{r}(M)]$ defines an inclusion of vector spaces.

Using the description of simple modules of equivariantization from \cite{buna} it follows that in fact  $K_0(\cc^F)= A^F$. Indeed, by \cite[Remark 3.12]{buna} a $k$-linear basis of $K_{0}(\cc^{F})$ is given by the elements  $\sum_{r \in F/F_{Y}}T^{r}(Y)$ where $Y$ runs throug all the orbits of the isomorphisms classes of simple objects of $\cc$. Recall that  $F_{Y}:=\{x \in F\;|T^{x}( Y)\cong Y\}$ is the inertia subgroup.

On the other hand if $M \in \cc_{g}^{F_{g}}$ is a simple object then by \cite[Theorem 2.12]{buna} it follows that $M\simeq \Ind^{F_{g}}_{F_{Y}}(Y \ot \pi)$, for a simple object $Y\in \cc$,  a constituent of $M$ and some projective representation $\pi$ of $(F_{g})_{Y}:=F_{Y}\cap F_{g}$. Since a $k$-linear bases of  $A^{F}$ is given by $\sum_{r\in F/F_{g}}T^{r}(M)$ where $M \in \cc_{g}^{F_{g}}$ is a simple object it follows from the above description of $M$ that the same vectors $\sum_{r \in F/F_{Y}}T^{r}(Y)$ form also a $k$-linear basis of $A^{F}$.

It remains to show that the multiplication from Theorem \ref{tpgr}  coincides to the multiplication described in Equation \eqref{multo}.  Once this is proven, it follows that condition $(\mathrm{C}5)$ is also satisfied since the multiplication in $K_{0}(\cc^{F})$ is associative.

Note that as in \cite{scoh} one has that $\co(g)\co(h)=\sqcup_{x \in F_{h }\slash F \backslash F_{g}}\co(\;^{x}gh)$. On the other hand multiplication formula from Theorem \ref{tpgr} shows
\beqarn
[S(g, M)][S(h, N)] & = & \sum_{x \in D}\sum_{r \in F/F_{\;^{x}gh}}[T^{r}(m_{^{x}g, h}(T^{x}(M), N))]\\ &=&  \sum_{x \in D}\sum_{r \in F/F_{\;^{x}gh}}m_{^{rx}g, \;^{r}h}([T^{rx}(M)], [T^{r}(N)])%\\ &=& \sum_{x \in D}\sum_{r \in F/F_{\;^{x}gh}}
\eeqarn
We have to show that this multiplication coincides to the one given in Equation \eqref{multo}.
For a fixed $x\in D$ note that $m_{^{rx}g, \;^{r}h}([T^{rx}(M)], [T^{r}(N)])$ with $r \in F/F_{\;^{x}gh}$ runs through all the orbit of of $\;^{x}gh$. Thus the term in $A(\;^{x}gh)$ of the above product coincides to $$\sum _{\{y \in D\:|\co(\;^{y}gh)=\co(\;^{x}gh)\}}m_{^{r_{y}y}g,\;^{r_{y}}h}([T^{r_{y}y}(M), T^{r_{y}}(N))$$ where $r_{y}\in F/F_{\;^{y}gh}$ is uniquely chosen such that $\;^{r_{y}y}g\;^{r_{y}}h=\;^{x}gh$.

One needs to show that the set $\{y \in D\:|\co(\;^{y}gh)=\co(\;^{x}gh)\}$ has the same cardinality as the set $[(\co(g)\times \co(h))\cap \{(a, b)\in G\times G\:|ab=\;^{x}gh\}]/F_{\;^{x}gh}$. In order to do this we construct a bijection between these two sets. If $\;^{r}g\;^{s}h=\;^{x}gh$ then send the orbit of the pair $(\;^{r}g,\;^{s}h)$ to the double coset $F_{h}s^{-1}rF_{g}$. Clearly this map is well defined. Conversely, define $F_{h}yF_{g}\mapsto \co((\:^{ry}g, \;^{r}h))$ where $r$ is chosen such that  $\:^{r}(\;^{y}gh)=\;^{x}gh$. It is easy to check that these two maps are one inverse to another.
\epf

\bl Suppose that $F$ acts coherently on a $G$-graded fusion category $\cc$ with respect to a given action by group automorphisms of $F$ on $G$.
Then with the above notations it follows that
\beq
T^{x}(m_{g, h}(M, N))\simeq m_{\;^{x}g, \;^{x}h}(T^{x}(M), T^{x}(N))
\eeq
for all $M\in \cc^{F_{g}}$ and $N\in \cc^{F_{h}}$.
\el
\bpf
As objects of $\cc$ one has that \beqn m_{g, h}(M, N)=\oplus_{r \in F_{gh}/F_{g}\cap F_{h}}T^{r}(M\ot N).\eeqn On the other hand \beqn m_{\;^{x}g, \;^{x}h}(T^{x}(M), T^{x}(N))=\oplus_{r \in F_{gh}/F_{g}\cap F_{h}}T^{xr^{-1}x}(T^{x}(M)\ot T^{x}(N)).\eeqn as objects of $\cc$.

It can be checked directly that $F: T^{x}(m_{g, h}(M, N))\ra m_{\;^{x}g, \;^{x}h}(T^{x}(M), T^{x}(N))$ given on components by
\beqn
T^{x}(T^{r}(M\ot N)) \xra{(T_{2}^{x,r})_{M\ot N}} T^{xr}(M\ot N)\xra{(T_{2}^{xrx^{-1},x})^{-1}_{M\ot N}} \eeqn \vskip -0,5cm \beqn \xra{(T_{2}^{xrx^{-1},x})^{-1}_{M\ot N}} T^{xrx^{-1}}(T^{x}(M\ot N))\xra{T^{xrx^{-1}}((T_{2}^{x})^{M, N})} T^{xrx^{-1}}(T^{x}(M)\ot T^{x}(N))
\eeqn
is an isomorphism in $\cc^{F_{\;^{x}(gh)}}$.

If $z' \in F_{gh}$ note that the equivariant structure of $m_{g, h}(M, N)$ is given on the components by 
\beqn
\mu^{z, r}_{m_{g, h}(M, N)}:T^{z}(T^{r}(M\ot N))\xra{(T_{2}^{z, r})_{M\ot N}} T^{zr}(M\ot N) \xra{(T_{2}^{r', h})^{-1}_{M\ot N}}\eeqn \vskip -0,5cm \beqn  T^{r'}(T^{l}(M\ot N))
\xra{T
^{r'}((T^{l}_{2})^{M, N})} T^{r'}(T^{l}(M)\ot T^{l}(N))\xra{T^{r'}(\mu_{M}^{l}\ot \mu_{N}^{l})} T^{r'}(M\ot N) \eeqn \vskip -0,5cm \beqn 
\eeqn
where $zr=r'l$ with $l \in F_{g}\cap F_{h}$ and $r' \in F_{gh}/F_{g}\cap F_{h}$.
Using Equation \eqref{conjmu} it follows that the equivariant structure of $T^{x}(m_{g, h}(M, N))$ is given on components by
\beqn
\;^{x}\mu^{xzx^{-1},\; r}_{T^{x}(m_{g, h}(M, N))}:T^{xzx^{-1}}(T^{x}(T^{r}(M\ot N))\xra{(T_{2}^{xzx^{-1}, x})_{M\ot N}} T^{xz}(T^{r}(M\ot N))\xra{} \eeqn \vskip -0,5cm \beqn \xra{(T_{2}^{x, z})^{-1}_{M\ot N}} T^{x}(T^{z}(T^{r}(M\ot N)) \xra{T^{x}(\mu^{z, r'}_{m_{g, h}(M, N)})}T^{x}(T^{r'}(M\ot N)).
\eeqn
On the other hand the equivariant structure $\nu^{xzx^{-1},\; r}_{ _{m_{\;^{x}g, \;^{x}h}(T^{x}(M), T^{x}(N))}}$ of the object $m_{\;^{x}g, \;^{x}h}(T^{x}(M), T^{x}(N))$
is given on the components by
\beqn
T^{xzx^{-1}}(T^{xr^{-1}x}(T^{x}(M)\ot T^{x}(N))) \xra{(T_{2}^{xzx^{-1}, xrx^{-1}})_{M\ot N}} T^{xzrx^{-1}}(T^{x}(M)\ot T^{x}(N))
\xra{} 
%\eeqn \vskip -0,5cm
\eeqn
\beqn \xra{(T_{2}^{x, z})^{-1}_{M\ot N}} T^{x}(T^{z}(T^{r}(M\ot N)) 
\xra{T^{x}(\mu^{z, r'}_{m_{g, h}(M, N)})}T^{x}(T^{r'}(M\ot N)).
\eeqn
Therefore verifying that $F$ is an isomorphism resumes to the commutativity of the following diagram (D1) made of solid arrows below. The commutativity of the bottom left rectangle of diagram (D1) follows from commutativity of diagram (D2). For shortness the maps in the diagrams are omitted but they are all uniquely determined from the group action of $G$ on $\cc$.
\subsubsection{Grothendieck rings of abelian cocentral extensions} In \cite[Theorem 4.8]{scoh} it is shown that the Grothendieck rings $\mathcal{G}_{0}(H)$ associated to abelian cocentral extensions have the multiplication structure given in Equation \eqref{multo}. Proposition \ref{cocsp} implies that the same result holds for any cocentral extension of semisimple Hopf algebras.

\br Note that compatibility condition $(\mathrm{C2})$ is not stated in \cite{scoh} on page 5 although it is stated as a property for Grothendieck groups of cocentral extensions on the last page of the paper. \er

\subsubsection{On the Grothendieck ring of the center of a fusion category} Suppose that $\cc$ is a $G$-graded fusion category $\cc=\oplus_{g \in G}\;\cc_{g}.$
Then by \cite[Theorem 4.1]{gnn} its Drinfeld center $\cz(\cc)\simeq \cz_{\cc_{1}}(\cc)^{G}$, the equivariantization of the relative center $\cz_{\cc_{1}}(\cc)$ by a certain action of the finite group $G$. Moreover \cite[Theorem 3.2]{gnn} shows that the relative center $\cz_{\cc_{1}}(\cc)$ is a $G$-crossed braided fusion category. In view of Example \ref{gbb} one can apply Theorem \ref{idr}. It follows that the Grothendieck ring of $\cz(\cc)$ has the ring structure described in Equation \eqref{multo}. Note that for the Grothendieck ring of a Drinfeld double of a semisimple Hopf algebra this description was already obtained in \cite{gdd}.

\begin{center}
\begin{equation}\tag{D1}
{\tiny
\begin{tikzpicture}
  \matrix (m) [matrix of math nodes, row sep=5. 5 em,column sep=1em,minimum width=2em]
  {
     T^{xzx\inv}(T^{x}(T^{r}(M\ot N))) & 
     T^{xzx\inv}(T^{xr}(M\ot N)) &   
     T^{xzx\inv}(T^{xrx^{-1}}(T^{x}(M\ot N))) \\%1
    T^{xz}(T^{r}(M\ot N)) & & T^{xzx\inv}(T^{xrx^{-1}}(T^{x}(M)\ot T^{x}(N))) \\%2
    T^{xzr} (M\ot N) & T^{xzrx^{-1}} (T^{x}(M\ot N))&  T^{xzrx^{-1}} (T^{x}(M)\ot T^{x}(N)) \\%3
     T^{x}(T^{zr} (M\ot N) )&  T^{xr'x^{-1}} (T^{xlx^{-1}}(T^{x}(M\ot N))  & T^{xr'x^{-1}} (T^{xlx^{-1}}(T^{x}(M)\ot T^{x}(N))) \\%4
     T^{x}(T^{r'}(T^{l}(M\ot N)))&  T^{xr'x^{-1}}(T^{xl}(M\ot N))  &   T^{xr'x^{-1}}(T^{xlx^{-1}}(T^{x}(M))\ot T^{xlx^{-1}}(T^{x}(N)))\\%5 
    T^{x}(T^{r'}(T^{l}(M)\ot T^{l}(N))) &  T^{xr'x^{-1}}(T^{x}(T^{l}(M\ot N))) & T^{xr'x^{-1}}(T^{xl}(M)\ot T^{xl} (N))   \\%6
     T^{x}(T^{r'}(M\ot N))) &   T^{xr'x^{-1}}(T^{x}(T^{l}(M)& T^{xr'x^{-1}}(T^{x}(T^{l}(M)\ot T^{x}(T^{l} (N))))\\%7 
     T^{xr'}(M\ot N)& T^{xr'x^{-1}}(T^{x}(M\ot N)) & T^{xr'x^{-1}}(T^{x}(M)\ot T^{x}(N)) .\\%8
     };
     \path[-stealth]
     (m-1-1) edge node [above] {$$}(m-1-2)%A1horiz T^{xzx^{-1}}((T^{x,r}_{2})_{M\ot N})
     (m-1-2) edge node [above] {$$}(m-1-3) %T^{xzx^{-1}}((T^{xrx^{-1},x}_{2})^{-1}_{M\ot N})
     (m-8-1) edge node [below] {$$}(m-8-2)%A8horiz %(T^{xr'x^{-1}, x}_{2})^{-1}_{M\ot N}
      (m-8-2) edge node [below] {$$}(m-8-3)% T^{xr'x^{-1}}((T^{x}_{2})^{M, N})
        (m-3-1) edge [dashed] node  [above] {$$}(m-3-2)%A8horiz (T^{xzrx^{-1}, x}_{2})_{M\ot N}
        (m-3-2) edge [dashed] node [above] {$$}(m-3-3)%T^{xzrx^{-1}}((T^{x}_{2})^{M, N})
       % (m-4-1) edge [dashed] node  [above] {hello}(m-4-2)%A8horiz
        (m-4-2) edge [dashed] node [above] {}(m-4-3)
       % (m-5-1) edge [dashed] node  [above] {}(m-5-2)%A8horiz
        (m-5-2) edge [dashed] node [above] {}(m-6-3)
        %(m-6-1) edge [dashed] node  [above] {}(m-6-2)%A8horiz
       % (m-7-3) edge [dashed] node  [above] {$\mu^{l}_{M\ot N}$}(m-7-1)%A8horiz
        (m-7-2) edge [dashed] node [above] {}(m-7-3)
      (m-8-2) edge  node [above] {}(m-8-3)
      %%%%% verticals\
       (m-1-1) edge node [left] {$$} (m-2-1)%vertical (T_{2}^{xzx^{-1}, x})_{T^{r}(M\ot N)}
       (m-2-1) edge node [left] {$$} (m-3-1)%vertical (T_{2}^{xz, r})_{M\ot N}
       (m-3-1) edge node [left] {$$} (m-4-1)% (T_{2}^{x, zr})^{-1}_{M\ot N}
       (m-4-1) edge node [left] {$$} (m-5-1)%T^{x}((T_{2}^{r', l})^{-1}_{M\ot N})
       (m-5-1) edge node [left] {$$} (m-6-1)% T^{x}(T^{r'}((T^{l}_{2})^{M, N}))
     (m-6-1) edge node [left] {$$} (m-7-1)% T^{x}(T^{r'}(\mu^{l}_{M} \ot \mu^{l}_{N})
     (m-7-1) edge node [left] {$$} (m-8-1)%vertical (T^{x, r'}_{2})_{M\ot N}
     (m-1-3) edge node [right] {$$} (m-2-3)%vertical T^{xzx^{-1}}(T^{xrx^{-1}}(T^{x}_{2})^{M, N})
       (m-2-3) edge node [right] {$$} (m-3-3)%vertical (T_{2}^{xzx^{-1}, xrx^{-1}})_{T^{x}(M)\ot T^{x}(N)}
       (m-3-3) edge node [right] {} (m-4-3)%
       (m-4-3) edge node [left] {$$} (m-5-3)%
       (m-5-3) edge node [left] {} (m-6-3)%{\tiny $T^{xr'x^{-1}}((T^{xlx^{-1}, x}_{2})^{-1}_{M}\ot (T^{xlx^{-1}, x}_{2})^{-1}_{N})$}
     (m-6-3) edge node [left] {$$} (m-7-3) % T^{xr'x^{-1}}((T^{x, l}_{2})^{-1}_{M}\ot (T^{x, l}_{2})^{-1}_{N})
     (m-7-3) edge node [left] {$$} (m-8-3) %T^{xr'x^{-1}}(T^{x}(\mu^{l}_{M} \ot \mu^{l}_{N}))
     (m-1-2) edge [dashed] node [right] {$$} (m-3-1)% (T^{xzx^{-1}, xr}_{2})_{M\ot N}
        (m-1-3) edge [dashed] node [left] {$$} (m-3-2)
        (m-3-2) edge [dashed] node [left] {$$} (m-4-2)
          (m-4-2) edge [dashed] node [left] {$$} (m-5-2)
          (m-5-2) edge [dashed] node [left] {$$} (m-6-2)
           (m-6-2) edge [dashed] node [left] {$$} (m-7-2)
            (m-7-2) edge [dashed] node [left] {$$} (m-8-2)
            %(m-6-2) edge [dashed] node [left] {$$} (m-7-2)
            ;
  \end{tikzpicture}
}
\end{equation}
\end{center}

\begin{center}
\begin{equation}\tag{D2}
\begin{tikzpicture}
  \matrix (m) [matrix of math nodes, row sep=4.5 em,column sep=1 em, minimum width=2em]
  {
    T^{xzr} (M\ot N) & & T^{xzrx^{-1}} (T^{x}(M\ot N)) \\%1
     T^{x}(T^{zr} (M\ot N) )& & T^{xr'x^{-1}} (T^{xlx^{-1}}(T^{x}(M\ot N))  \\%2
     T^{x}(T^{r'}(T^{l}(M\ot N)))&  T^{xr'}(T^{l}(M\ot N)) & T^{xr'x^{-1}}(T^{xl}(M\ot N)) \\%3 
    T^{x}(T^{r'}(T^{l}(M)\ot T^{l}(N))) & T^{xr'}(T^{l}(M)\ot T^{l}(N) )& T^{xr'x^{-1}}(T^{x}(T^{l}(M\ot N)))   \\%4
     T^{x}(T^{r'}(M\ot N))) & &  T^{xr'x^{-1}}(T^{x}(T^{l}(M)\\%5 
     T^{xr'}(M\ot N)& & T^{xr'x^{-1}}(T^{x}(M\ot N)) \\%6
     };
      \path[-stealth]
     (m-1-1) edge node [left] {$$} (m-2-1)%vertical
       (m-2-1) edge node [left] {$$} (m-3-1)%vertical
       (m-3-1) edge node [left] {$$} (m-4-1)%
       (m-4-1) edge node [left] {$$} (m-5-1)%
       (m-5-1) edge node [left] {$$} (m-6-1)%
        (m-1-3) edge node [left] {$$} (m-2-3)%vertical
       (m-2-3) edge node [left] {$$} (m-3-3)%vertical
       (m-3-3) edge node [left] {$$} (m-4-3)%
       (m-4-3) edge node [left] {$$} (m-5-3)%
       (m-5-3) edge node [left] {$$} (m-6-3)%
        (m-1-1) edge node [above] {$$}(m-1-3)%A1horiz%A6 horiz
        %dashed inside
        (m-6-1)edge node [above] {$$}(m-6-3)
         (m-1-1) edge [dashed] node [above] {$$}(m-3-3)
         (m-1-1) edge [dashed] node [above] {$$}(m-3-2)
         (m-3-1) edge [dashed] node [above] {$$}(m-3-2)
         (m-3-2) edge [dashed] node [above] {$$}(m-4-3)
         (m-4-1) edge [dashed] node [above] {$$}(m-4-2)
         (m-3-2) edge [dashed] node [above] {$$}(m-4-2) 
         (m-4-2) edge [dashed] node [above] {$$}(m-5-3)
         (m-4-2) edge [dashed] node [above] {$$}(m-6-1)           ;
  \end{tikzpicture}
  \end{equation}
  \end{center}
\epf
%Indeed, if $\cc=\oplus_{g \in G}\cc_{g}$ with $\cc_{1}=\cd$ then using the fact that $\cz(\cc)\simeq\cz_{\cd}(\cc)^G$ and the category $\cz_{\cd}(\cc)$ is a braided $G$-crossed category by \cite{gnn}.
\bibliographystyle{amsplain}
\bibliography{green-funct}
\ed
\section{Indecomposable equivariantized objects of tensor categories} Suppose that $G$ acts linearly on the abelian category $\cc$.

\subsection{Group actions on Hom-spaces}
Let $M, N \in \cc^{G}$ be two equivariantized objects. There is an action of G on $\Hom_{\cc}(M, N)$ given as follows
\bq
g.f=\mu^{g}_{N}(T^{g}(f))(\mu^{g}_{M})^{-1}
\eq
for all $f \in \Hom_{\cc}(M, N)$ and $g \in G$. It is easy to check that
\bq
g.(uv)=(g.u)(g.v)
\eq
for any $M, N, P\in \cc^{G}$ and any two morphisms $M\xra{v} N\xra{u} P$ in $\cc$.
\md
\bl Let $M, N \in \cc^{G}$. Then
\beq
 \Hom_{\cc}(M, N)^{G}=\Hom_{\cc^{G}}(M, N)
 \eeq
\el
\subsubsection{Properties of induction and restriction functors} Let $H\leq G$ be a subgroup and $W\in \cc^{G}$. There is an epimorphism in $\cc^{H}$
\bq\label{indr}
\ind^{G}_{H}(\res^{G}_{H}(W)) \xra{\al_{W}} W 
\eq
given by the  canonical projection $\oplus_{t \in G/H}T^{t}(W)\ra W$. It is easy to verify that this is a morphism in $\cc^{H}$. Indeed, for any $g \in G$ the equivariant structure of $\ind^{G}_{H}(\res^{G}_{H}(W))$ on the component $\ro^{1}(W)$ is given by:
\bq
T^{g}(T^{1}(W))=T^{g}(W)\xra{(T^{r,h}_{2})^{-1}} T^{r}(T^{h}(W))\xra{} T^{r}(W)
\eq
\subsubsection{}
Moreover for any $U \in \cc^{H}$ there is also a canonical projection in $\cc^{H}$
\bq\label{retri}
\res^{G}_{H}(\ind^{G}_{H}(U)) \xra{\beta_{U}} U 
\eq
\blue{given by the canonical projection}
\subsubsection{The endomorphism sum in $\mtr{End}_{\cc^{G}}(U)$}
Note that there is an endomorphism of $\ind(U)$ such that
\bq
\sum_{i}x_if=1_{U}
\eq
for all $f \in \End(U)$.
\subsubsection{Mashcke theorem for equivariantized tensor categories.}

\bt Let $G$ be a finite group acting on a tensor category $\cc$. Then the tensor category $\cc^G$ is a fusion category if and only if the characteristic of $k$ does not divide the order of $G$ and the category $\cc$ is a semisimple category.
\et

\bpf Suppose that $\cc$ and $\rep(G)$ are fusion categories.
Let $\pi: M \ra N$ be an epimorphism in $\cc^{G}$ and $h$ be a splitting of $\pi$ in $\cc$. Then define $h' : N \ra M$ by
\bq
h'= \frac{1}{|G|}\sum_{g \in G} g.h
\eq

Then $h'$ is a morphism of G equivariant fusion categs since $g.h'=h'$ for all $g\in G$. Moreover $h'$ splits $\pi$ since
\bq
\pi h'=\frac{1}{|G|}\sum_{g\in G} \pi g.h=\frac{1}{|G|}\sum_{g \in G} g(\pi h)=\frac{1}{|G|}
\sum_{g \in G} g. 1_{N}=1_{N}
\eq

The converse also follows. Since $F: \cc^G \ra \cc$ is a surjective tensor functor one has that $\cc$ is semisimple. Moreover $\rep(G)$ is a tensor subcategory of $\cc$ and thus semisimple. Then classical Mashke's theorem implies that  order of $G$ does does not divide the characteristic of $k$.
\epf

\subsection{Vertices and sources}
\blue{Definition of a split sequence}
\subsubsection{Relatively $H$-free objects}
 is relatively $H$-free if for any object 
\subsubsection{Relatively $H$-projective objects}
%\bn{def}
Let $(W, \mu)\in \cc^{G}$ and $H$ a subgroup of $G$. Then $W$ is
called $H$-projective if, whenever there is a short exact sequence in $\cc^{G}$ of the form 
\beq
0\ra V\ra U\ra W\ra 0
\eeq
such that under restriction to $\cc^{H}$ the exact sequence splits then it also splits in $\cc^{G}$.
%\end{def}
\md
We generalize the following characterization of $H$-projectivity which  is due to Higman in the case of representations of finite groups.
\bt Let H be a subgroup of the group G. Then each of the following properties for any $ W\in \cc^{G}$ are equivalent:
\md
(a) $W$ is H-projective;
\md
(b) $W|\ind^{G}_{H}(W)$
\md
(c) $W|\ind^{G}_{H}(U)$ for some object $U\in \cc^{H}$ 
\md
(d) there is $f\in \mtr{End}_{\cc^{H}}(W)$ such that $\sum_{i}x_{i}.f=1_{W}$.
\et
\bpf
$a)\implies b)$ There is a surjective morphism $\ind^{G}_{H}(W)\ra W$ define as... that splits when restricted to $\cc^{H}$.
\md Trivially $b)\implies c)$.
\md
$c)\implies a)$
\epf
%(d) there exists f~ End,,( W) such that C;= , x; ' f x i = 1,
if possible to write Mackey for non semisimple fusion categories

\subsection{Vertices and sources}
Let $V \in \cc^{G}$ be an indecomposable $G$-equivariant object. Then a subgroup $Q$ of $G$ is a
vertex of $V$ if $V$ is $Q$-projective but $V$ is not $H$-projective for any proper subgroup of $Q$. 
\md

It follows from Theorem 5.l(c) that corresponding to any vertex
$Q$ of $V$ there is at least one indecomposable $Q$-equivariant object $U$ such that $V| U^{G}$ ;
such an Q-equivariant object is called a source of $V$.
%%%
\bc Let $G$ acting on a tensor category $\cc$ and $H,L$ be subgroups of the group G.
(i) If $x^{{-1}}Hx\subset  L$ for some $x \in G$, then each $H$-projective $G$-equivariant object is $L$-projective.
(ii) IF $H \subset $L and $V$ is an H-projective $L$-equivariant object, and W is an AC-equivariant object such that $W | V^{G}$, then $W$ is an $H$-projective $G$-equivariant object.
\ec
\bpf(i) Let $W$ be an $H$-projective $G$-equivariant object. By part (c) of the
theorem there exists an AH-equivariant object U such that W I UG.T hen U �3 x is an
A(x- 'Hx)-subequivariant object of U" and it is easily seen that (U �3 x)" = UG as
$G$-equivariant objects. Let V be the induced AL equivariant object (U @I x)". Since inducing is a
transitive operation, p = (U �3 x)" = UG and so W I p. Thus the theorem
shows that W is Lprojective.
(ii) Since V is an H-projective AL-equivariant object, part (c) of the theorem shows
that $V | U^{G}$ for some AH-equivariant object $U$. Then $V|  ( U^{L})^{G} = U^{G}$, so $W | U^{G} $ and hence $W$ is $H$-projective.
\epf
\md
\bl
Let $H$ be a normal subgroup of index $p$ in the group $G$, and
suppose that $k$ is a perfect field of characteristic $p$. Let $W$ be an absolutely indecomposable object of $\cc^{H}$, and suppose that the indecomposable components
of $\ind^{G}_{H}(\res^{G}_{H}(H))$ are all isomorphic to $W$. Then $\ind^{G}_{H}(W)$ is an absolutely indecomposable object of $\cc^{G}$.
\el
\bpf
Let $R$ be the $k$-algebra  $\mtr{End}_{\cc}(W)$.
\md
One can define an algebra isomorphism
\bq
\psi:\End_{\cc^G}(\ind^G_H(W))\ra \text{Mat}(p, R)
\eq
defined by $a \mapsto [\al_{ij}]$  where $\al_{ij} \in \End_{\cc}(W)$ are defined by $\al_{ij}:=\psi^i a_{ij}\psi^{-j}$. Here the morphism $a_{ij}$ is the $(i,j)$-component of $a$ defined $T^{y^i}(W)\ra T^{y^j}(W)$.
\md
Note that any element $g \in G$ define s a morphism of $L_g \in \End_{\cc}(\ind^G_H(W))$ defined by
{\small
\bq
T^g(T^{y^i})(W))\xra{} T^{ga}(W) \xra{} T^{y^jh}(W) \xra{} T^{y^j}(T^h(W))\xra{(T^{y^j}(\mu^W_h
)} T^{y^j}(W) \xra{\phi^i} W\eq}

It can be shown that
$\psi(h)=$ and $\psi(y)=$...
\epf
\blue{In Curtis and Reiner there is a different proof using graded rings.}
\md
 With this Lemma one can prove an analogue of Green's theorem:
\bt
 Let $H$ be a normal subgroup of index $p$ in a group $G$. Let $k$ be
a field of characteristic $p$ , and let $W$ be an absolutely indecomposable object of $\cc^{H}$.
Then $\ind^{G}_{H}(W)$ is an absolutely indecomposable object in $\cc^{G}$.
\et
\bpf
It is enough to prove the theorem for an algebraically closed field of characteristic $p$.
\md
Consider the set $Z_{W}:=\{x \in G\:|T^{x}(W)\simeq W \;\text{in}\;\cc^{H}\}$ This is a subgroup of $G$ containing $H$ and therefore we have the following two cases:
\md{\it Case 1:} Suppose that $Z_{W}=H$ and let $y\in G\slash H$. Then
\bq
\res^{G}_{H}(\ind^{G}_{H}(W))=\oplus_{i=0}^{{p-1}}T^{y^{i}}(W)
\eq
Using the Krull-Schmidt Theorem we know that $\ind^{G}_{H}(W)$ has an indecomposable component $V$ such that $W|\res^{G}_{H}(V)$ in $\
\cc^{H}$. But then $T^{y^{i}}(W)| T^{y^i}(V)$ and the latter equals $V$, since $V \in \cc^G$. Since $W$ is indecomposable, the objects $T^{y^i}(W)$ with $i=0, \cdots , p - 1$ are indecomposable, and they are mutually nonisomorphic because
$I( W) = H$. Therefore the Krull-Schmidt theorem shows that the sum
of these modules divides $V$; that is, $\ind^G_{H}(W )| V$. But this implies $\ind^G_{H}(W )=V$,  and so
$\ind^g(H(W))$ is indecomposable
\md {\it Case 2:} $Z_{W}=G$. In this case we are in the situation of the previous Lemma and the result follows from there.
\epf
\blue{You can also ask Kleiner if there is Green's Theorem for Green functors}
\bibliographystyle{amsplain}
\bibliography{bob}
\ed

\bibliographystyle{amsplain}
\bibliography{bob}
\ed
\section{coherent actions of groups on graded  fusion categories}\lb{hmg} For any category $\cc$ denote by $\co(\cc)$ the set of all objects of $\cc$. Suppose that $\cc$ is a fusion category graded by a finite group $G$. Recall that this means that
\beq
\cc=\oplus_{g \in G}\cc_{g}
\eeq
as abelian categories, and the tensor functor $\ot: \cc \times \cc\ra \cc $ sends $\cc_g\ot \cc_h$ into $\cc_{gh}$. For an object $V \in \cc$ define by $V_g$ the homogenous component of degree $g$ from the above grading of $\cc$.
 \md
 
Suppose further that another finite group $F$ acts by group automorphisms on $G$. Suppose that $F$ also acts by tensor automorphisms on the category $\cc$ via the action $T:F\ra \underline{\mtr{Aut}}_{\ot}(\cc)$ given by $x \mto T^x  :\cc \ra \cc$.

\bn{defn} An action of a finite group $F$ on the $G$-graded fusion category $\cc$ is called {\it coherent}
 if 
\beq
T^x(\cc_g)\subset \cc_{\;^xg}
\eeq
for all $x \in F$ and $g \in G$.
\end{defn}
\md
\br
Note that the action of $F$ is coherent if and only if $\cc_1$ is stable under the action of $F$. Indeed let $V, W\in \cc_g$ and suppose that $T^x(V)\in \cc_{g_1}$ and $T^x(W) \in \cc_{g_2}$. Then $T^x(V^*\ot W)\in \cc_1$ since $V^*\ot W \in \cc_1$. On the other hand $T^x(V^*\ot W)\simeq T^x(V^*)\ot T^x(W)\in \cc_{{g_1}^{-1}g_2}$ which implies that $g_1=g_2$. Denoting $\;^xg:=g_1$ then it is easy to check that this defines an action of $F$ on $G$ by group automorphisms.
\er
\subsection{Examples of coherent actions of groups and their equivariantized categories}
\bn{example}\lb{gbb}{\it Braided $G$-crossed categories.}\end{example}\md
Recall that a crossed fusion category is a quadruple $(\cc, G, T, c)$, where $G$ is a finite group, $\cc$ is a tensor category with a (not necessarily faithful) $G$-grading $\cc=\oplus_{{g \in G}}\cc_{g}$
 and a tensor action $T : G \ra \underline{Aut}_{\ot}(\cc)$, $g \mapsto T_{g}$ satisfying $T_{g}(\cc_{h})\subseteq \cc_{ghg\inv}$. Moreover the crossed braiding $c$ defined by $c(X, Y):X\ot Y \ra T_{g}(Y)\ot X$ for all $X\in \cc_{g}$ and $Y \in \cc$. The compatibility conditions that have to be satisfied this quadruple can be found for example in \cite{DGNO}.\md Suppose that $\cc$ is a braided $G$-crossed category. We see that in this case $F=G$ acts homogeneously on $\cc$ where the action of $F$ on $G$ is given by conjugation.
\bn{example}{\it 
Cocentral extensions of semisimple Hopf algebras}\md
\end{example}
%\red{One can also put it as an appendix}
Suppose that we have a cocentral extension of Hopf algebras 
\beq\lb{cocen}
k \ra B\xrightarrow{i} H\xrightarrow{\pi} kF\ra k
\eeq

Recall that this means the $kF^*\subset \mtc{Z}(H^*)$ via $\pi^*$. Following Proposition 3.5 of \cite{natalecoc} it follows that for any such extension one has that $F$ acts on $\Rep(B)$ and $\Rep(H)=\Rep(B)^F$.

On the other hand, using the reconstruction theorem from \cite{AD} it follows that 
\beq\lb{andc}
H \simeq B\;^{\tau}\#_{\sg} \;kF
\eeq
for some cocycle $\sg:B\ot B \ra kF$ and some dual cocyle $\tau:kF \ra B\ot B$. 
\md

Recall that in this situation the multiplication on $H$ is given by 
%\beq(a\#_{\sg}h)(b\#_{\sg}l)=a(h.b)\#hl\eeq
\bq\lb{mu}
(b\#_{\sg}f)(c\#_{\sg}g)=b(f.c)\sg(f,g)\#_{\sg}fg
\eq
and the comultiplication on $H$ is given by:
%\beq\Delta(b\#g)=(b_1\tau(g)_j\#_{\sg}g)\ot (b_2\tau(g)^j\#_{\sg}g)\eeq
\bq\lb{com}
\D(b\#_{\sg}\bar{f})=(b_1\tau(\bar{f})_i\#_{\sg}\bar{f})\ot(b_2\tau(\bar{f})^i\#_{\sg}\bar{f})
\eq
Moreover the antipode is given by the formula $$S(a\#_{\sg}g)=g^{-1}S(a)=(g^{-1}.S(a))\#_{\sg}g^{-1}.$$
Note that here the weak action $kF\ot B \ra B$ is denoted by $f \ot b\mapsto f.b$ and we used the same notation as in \cite{AD}, 
\bq
\tau(f)=\tau(f)_i\ot \tau(f)^i\in B \ot B\in B \ot B.
\eq 
%and the weak cocation $T$ is denoted by $T(b)=\sum_i b_i\ot b^i \in B \ot kF$.

One may also suppose that for all $f \in F$ one has $\sg(f,f^{-1})=1$ and therefore $\bar{f}^{-1}=\bar{f^{-1}}$.  In particular, the cocycles $\sg: F\ot F \ra U(B)$ and $\tau:F \ra B\ot B$ should verify the following compatibility conditions:  

\bn{enumerate}
\item $\sg(1,g)=\sg(g,1)=1$
\item $(\eps_A\ot \id )(\tau(g))=1=(\id \ot \eps_A)(\tau(g))$
\item  \lb{cocyl} $[h.\sg(l,m)]\sg(h,lm)=\sg(h,l)\sg(hl,m)$
\item  \lb{twm} $h.(l.a)=\sg(h,l)[(hl).a]\sg^{-1}(h,l)$
\item $(\tau(g)_j)_1\tau(g)_l\ot (\tau(g)_j)_2\tau(g)^l\ot \tau(g)^j=\tau(g)_j\ot (\tau(g)^j)_1\tau(g)_r\ot  (\tau(g)^j)_2\tau(g)^r$

This means that $\tau(g)$ is a twist for $A^{op}$ or equivalently $\tau(g)^{-1}$ is a twist for $A$.
\item \lb{A} $\Delta_A(g.a) = \tau (g)(g.a_1 \ot g.a_2)\tau(g)^{-1}$ \\ (This is condition $(A)$ from \cite{AD})
\item \lb{D} $\Delta(\sg(g,h))\tau(gh)=\tau(g)(g.\tau(h)_p\ot g.\tau(h)^p)(\sg(g,h)\ot \sg(g,h))$
(This is condition $(D)$ from \cite{AD})  \end{enumerate} 
Note that conditions $(B)$ and $(C)$ from \cite{AD} are automatically satisfied for cocentral extensions.
%\subsubsection{Universal gradings for Hopf algebras}
\bp
With the above notations suppose that $\cc:=\rep(B)$ and let
\begin{equation}\label{univ}
\mtc{C}=\oplus_{g \in G}\mtc{C}_g
\end{equation}
be the universal grading of $\rep(B)$ where $K(B)=kG^*$ is the largest central Hopf subalgebra of $B$. Then $F$ acts on $G$ and the action of $F$ on $\cc$ is coherent with respect to the universal grading of $\cc$.
\ep
\bpf
First it will be shown that $F$ acts on the group $G$ by group automorphisms. In order to do that we first show that $F$ acts on $K(B)=kG^*$ via Hopf automorphism.
Note that if $b \in K(B)$ the $f.b \in K(B)$. Clearly $f.b \in \cz(B)$ since $F$ acts by algebra automorhisms on $B$.
\md
On the other hand using formula $(A)$ from \cite{AD} it follows that 
\beq
\Delta(f.b)=\tau(f)(f.b_1\ot f.b_2)\tau(f)^{-1}.
\eeq
Therefore if $b \in \cz(B)$ then 
$
\Delta(f.b)=(f.b_1\ot f.b_2)
$
which shows that $F. K(B)$ is a central Hopf subalgebra of $B$.  Thus $F.K(B)\subseteq K(B)$ and $F$ acts by Hopf algebra automorphisms on $K(B)$. Therefore there is an action of $F$ on $G$ such that the action of $F$ on $K(B)=kG^*$ is induced by
\bq
x.p_g=p_{\;^xg}
\eq
for all $x \in F$ and $g \in G$. Now to verify  that the action of $F$ on $\rep(B)$ is coherent one has to verify that for any $M \in \co(\cc_g)$ then $T^x(M)\in \co( \cc_{\;^xg})$. Recall from Proposition 3.5 of \cite{natalecoc} that $T^x(M)=M$ as vector spaces and the action of $B$ on $T^x(M)$ is given by
\bq
b.\;^xm=(x^{-1}.b)m
\eq
Thus for any $h \in G$ one has that 
\bq
{p_h.}\;^xm=(x^{-1}.p_h)m=p_{\;^{x^{-1}}h}\;m=\delta_{\;^{x^{-1}}h, g}\;m=\delta_{h, \;^{x}g}\;m
\eq
which shows that $T^x(M)\in \co(\cc_{\;^xg})$.
\epf
\bn{example}\end{example}\lb{twisted} 
Group actions on pointed fusion categories are always coherent. See \cite[Section 7]{naidu}. In particular,  it follows that the representation category of a (twisted) quantum double of a finite group is the equivariantization of a coherent action.

\sub{Simple objects for equivariantizations of coherent actions}With the 
settings from the beginning of this section note that the stabilizer $F_{g}$ of an element $g \in G$ acts by $k$-linear automorphisms on the abelian subcategory $\cc_{g}$ of $\cc$. In particular one obtains in this way an action by tensor automorphisms of  $F$ on the tensor category $\cc_{1}$.

\bl With the above notations, suppose that $(V, {\mu_x^V}_{\{x \in F\}})\in \cc^F$ is an equivariantized object with a canonical decomposition $V=\oplus_{g \in G}V_{g}$. Then for all $g \in G$ one has that $(V_g,  {\mu^V_x}_{\{x \in F_g\}}) \in \cc^{F_g}$.
\el

\bpf
Note that $\mu^V_x$ sends a component  $V_g$ to the component $V_{\;^{x^{-1}g}}$ since $\mu_x^V:T^x(V)\ra V$ is an isomorphism in $\cc$. Thus if $x \in F_g$ then $V_g$ is sent to itself by $\mu^V_x$. The compatibility conditions for an equivariantized object follow from the equivariancy of $V$.
\epf

Let $\Gm$ be a set of representatives for the orbits of the action of $F$ on the group $G$. For any $g \in G$ let $\co(g)$ be the orbit of $g$ under the action of $F$. Next proposition is a generalization of Proposition 2.7 of \cite{gnn}.

\bp Suppose that $F$ acts homogeneously on a graded fusion category $\cc$. Then simple objects of $\cc^F$ are parametrized by pairs $(g , M)$ where $g \in \Gm$ and $M \in {\cc_g}^{F_g}$  is a simple object.The simple object corresponding to the pair $(g, M)$ coincides to the induced object $\Ind_{F_{g}}^{F}(M)$.
\ep 

\bpf
Suppose that $(V, {\mu^{V}_{x}}_{\{x \in F\}})\in \cc^F$ is a simple object and let $V=\oplus_{g \in G}V_g$ with $V_{g }\in \cc_{g}$. If $V_g \neq 0$ then $T^x(V_g)=V_{\;^{x}g}$ is also not zero and therefore $\opl_{h \in \co(g)}V_h$ is an equivariantized object of $\cc^F$. Since $V$ is a simple object it follows that $V$ is supported only on the orbit $\co(g)$ of $g$, i.e. $V=\oplus_{{h \in \co(g)}}V_{h}$. Moreover it follows that $V_g\in \cc_g^{F_g}$ and one can associate to the object $(V, {\mu^{V}_{x}}_{\{x \in F\}})\in \cc^F$ the pair $(g, M)$ with $M:=F_g$. Clearly $M$ is simple if $V$ is simple.
\md
It remains to show that $V\simeq \ind^{F}_{F_{g}}(M)$. 
This is equivalent with showing that
\bq
V=\oplus_{r \in F/F_{g}}T^{r}(M)
\eq
Since $\mu_{x}^{V}:T^{x}(V) \xra{\simeq} V$ it follows that $\mu_{x}^{V}|_{T^{x}(V_{g})}:T^{x}(V_{g}) \xra{\simeq}V_{\;^{xg}}$ are also isomorphisms. This shows that $f: \ind_{F_{g}}^{F}(V_{g})\xra{\simeq} V$ is an isomoprhism in $\cc^{F}$ where $f:=\oplus_{r \in F/F_{g}}\mu_{x}^{V}|_{T^{x}(V_{g})}$.
\md
The converse also holds. If $M$ is a simple object of $\cc^{F_{g}}$ then $\ind_{{F_{g}}}^{F}(M)$ is simple. This can be seen simply just by inspecting the homogenous components of the induced module.
\epf
%}
\bc Suppose that $F$ acts homogeneously on a graded fusion category $\cc$ and
let $V=\oplus_{h \in \co(g)}V_{h}$ be a simple
object of $\cc^{F}$. Then 
\bq
\Ind^{F}_{F_{g}}(V_{g})
\simeq \ind_{F_{h}}^{F}(V_{h})
\eq
for all $h \in \co(g)$.
\ec
\bpf
Both terms are isomorphic to $V$.
\epf
\md
\subsection{Definition of $S(g, M)$} Let $M \in \cc^{F_{g}}$ be a simple object.
Define by $S(g, M)$ the simple induced object $\ind^{F}_{F_{g}}(M) \in \cc^F$. Thus as objects of $\cc$ one has that $S(g, M):=\oplus_{r \in F/F_g}T^r(M)$ with the equivariant $G$-structure from (\ref{muind}).

\br\lb{shift}\er Note also that if $(M , u^M) \in \cc_g^{F_g}$ then $(T^x(M), \;^xu^M)\in \cc_{\;^xg}^{F_{\;^xg}}$ with the equivariant structure $ \;^xu^M$ given by:
\beq
T^{xhx^{-1}}(T^x(M))\xra{T_2^{xhx^{-1}, x}} T^{xh}(M)\xra{(T_2^{x,h})^{-1}}T^x(T^h(M))\xra{T^x(u^M_h)} T^x(M)
\eeq

Therefore previous Corollary implies that:
\beq
S(g, M)\simeq S(\;^xg, T^x(M))
\eeq
for any $x \in F$, $g \in G$ and $M \in \cc_g^{F_g}$.
\md
%\subsection{Mackey's decomposition of an induced object}
\subsection{Adjoint tensor product formula}
\bp\lb{atp}Suppose that $\cc, \cd$ and $ \ce$ are tensor categories and $F_1:\cc \ra \cd$ and $F_2:\cc \ra \ce$ are two tensor functors with left adjoint functors $I_1:\cd\ra \cc$ and respectively $I_2:\ce\ra \cc$.  Then for any objects $M \in \co(\cd)$ and $N \in \co(\ce)$ one has a canonical isomorphism in $\cc$
\beq\lb{ti}
I_1(M)\ot I_2(N)\simeq I_2(F_2(I_1(M))\ot N).
\eeq

\ep
\bpf
It can be shown by a straightforward computation that 
\beq\lb{ti}
\Hom_{\cc}(I_1(M)\ot I_2(N), P)\simeq \Hom_{\cc}(I_2(R_2(I_1(M))\ot N)
, P)
\eeq
for any object $P \in \cc)$. This implies the conclusion.\epf
In particular one obtains that 
\beq\lb{ti2}
I_{1}(M)\ot V\simeq I_{1}(M\ot F_{1}(V))
\eeq for any objects $M\in \cd)$ and $V\in \cc)$.
\subsection{Tensor product formula}
Let $M\in \cc_g^{F_{g}}$ and $N \in \cc_h^{F_{h}}$ be two equivariant objects. Define the following equivariant object
\beq
m_{g,h}(M, N):=\ind_{F_{g}\cap F_{h}}^{F_{gh}}(\res^{F_{g}}_{F_{g}\cap F_h}(M)\ot \res_{F_{g}\cap F_h}^{F_h}(N))\in \cc_{gh}^{F_{gh}}.
\eeq
%\ind_{F_{\;^xg}\cap F_{h}}^{F_{\;^xgh}}(\res^{F_{\;^xg}}_{F_{\;^xg}\cap F_h}(\;^xM)\ot \res_{F_{\;^xg}\cap F_h}^{F_h}(N))
\bt \lb{tpgr} Suppose that $F$ acts homogeneously on the $G$-graded tensor category $\cc$. With the above notations one has
\beq\lb{tepg}
S({g, M})\ot S({h,N})\simeq \opl_{x \in D}S(\;^xgh, \;\;m_{ _{\;^xg, h}}(T^x(M),N))
\eeq
where $D$ is a set of representatives for the double cosets $F_g \backslash F / F_h$.
\et
\bpf
One has that $S(g, M)=\ind_{F_g}^F(M)$. Applying formula (\ref{ti}) one has that
\beqarn
S({g, M})\ot S({h,N}) & = & \ind_{F_g}^F(M)\ot \ind_{F_h}^F(N)\\
& \simeq & \ind_{F_h}^F(\res_{F_h}^F(\ind_{F_g}^F(M)\ot N)).
\eeqarn

On the other hand applying Theorem \ref{macky} one has that 
\beqarn
\res_{F_h}^F(\ind_{F_g}^F(M))\simeq \oplus_{x \in D}\ind_{{F_h} \cap \;^x{F_g}}^{F_h}(\res^{\;^x{F_g}}_{\;^x{F_g}\cap {F_h}}(T^x(M)))
\eeqarn
Then applying formula (\ref{ti2}) one obtains that
\beqarn
\res_{F_h}^F(\ind_{F_g}^F(M))\ot N\simeq \oplus_{x \in D}(\ind_{{F_h} \cap \;^x{F_g}}^{F_h}(\res^{\;^x{F_g}}_{\;^x{F_g}\cap {F_h}}(T^x(M))\ot N) \\ \simeq  \oplus_{x \in D}(\ind_{{F_h} \cap \;^x{F_g}}^{F_h}(\res^{\;^x{F_g}}_{\;^x{F_g}\cap {F_h}}(T^x(M)\ot \res^{{F_h}}_{\;^x{F_g}\cap {F_h}}(N))
\eeqarn
and therefore 
\beq
\ind_{F_g}^F(M)\ot \ind_{F_h}^F(N)\simeq  \oplus_{x \in D}\ind_{\;^{x}F_{g}\cap F_{h}}^{F}(\res^{\;^{x}F_{g}}_{\;^{x}F_{g}\cap F_h}(T^{x}(M))\ot \res_{\;^{x}F_{g}\cap F_h}^{F_h}(N))
\eeq
which by definition coincides to $\opl_{x \in D}S(\;^xgh, \;\;m_{ _{\;^xg, h}}(T^x(M),N))$.
\epf
\bc Suppose that $F$ acts homogeneously on the fusion category $\cc$. With the above notations it follows that
\beq\lb{ditf}
\cc^F \simeq \oplus_{g \in \Gm}\cc_g^{F_g}
\eeq
as indecomposable $\cc_{1}^F$-bimodule categories, where $\Gm$ is a set of representative elements for the orbits of the action of $F$ on $G$ .\ec
\bpf
Remark that $\cc_{1}^{F}$ is a tensor subcategory of $\cc^{F}$ consisting on those objects supported only on $\cc_{1}$. Define the functor $F: \cc^{F} \ra \oplus_{g \in \Gm}\cc_g^{F_g}$ by sending $S(g, M)\mapsto M$. Note that formula (\ref{tepg}) implies that 
\beqn
S(1, M)\ot S(h, N)\simeq S(h, \res^{F}_{F_{h}}(M)\ot N)
\eeqn
and
\beqn
S(h, N)\ot S(1, M)\simeq S(h, N\ot \res^{F}_{F_{h}}(M))
\eeqn
which shows that each $\cc^{F_{h}}$ is a $\cc_{1}^{F}$-bimodule category. The above formulae also show that $\cc^{F_{h}}$ is an indecomposable $\cc_{1}^{F}$-bimodule category for all $h \in H$.
\epf
\section{On the Grothendieck ring of a graded equivariantization}\lb{withsp}
\subsection{On the Green rings introduced by Witherspoon and Bouc}
%%%%%%%%%%%%%%%%
Let $F$ be a finite group acting by group automoprhisms on another finite group $G$ . Suppose  moreover that 
$
A=\oplus_{g \in G}A(g)
$
is a graded vector space endowed with two linear structures:
$
m_{g,h}:A(g)\ot A(h)\ra A(gh)
$
and 
$
c_{x,g}:A(g)\ra A(\;^xg).
$
Let $c_x:=\oplus_{g \in G}\;c_{x,g}:A\ra A$ and let
\beq\lb{inf}
A^F:=\{a \in A\;|\; c_x(a)=a\;\text{for all}\;x\in F\}
\eeq
be the subspace of $F$-invariants of $A$.
As in \cite{scoh} we suppose the following compatibility axioms are satisfied by the additional structures.
\beq\tag{H1}
c_1=\id \text{and\;} c_{xy}=c_x c_y.
\eeq 
\beq\tag{H2}
c_xm_{g,h}=m_{\;^xg,\;^xh}(c_x\times c_x)
\eeq
There is an element $1 \in A(1)$ such that
\beq\tag{H3}
c_x(1)=1 \; \text{for all\;} x \in L \; \text{and\;} m_{1,g}(1, \al_g)=m_{g,1}(\al_g,1)=\al_g
\eeq
for all $g \in G$. For any $g \in G$ let $F_g\subset F$ be the stabilizer subgroup of $g$. Then from Equation (H1) it follows that 
$F_g$ acts on $A(g)$. 
\beq\lb{multo}
(\al\beta)_g=\sum_{\{(h,k)\in L_g/G \times G\;| \;hk=g\}}m_{h,k}(\al_h,\beta_k)
\eeq
where the action of $F_g$ on $G\times G$ is diagonal.
Then it is shown in \cite{scoh} that the above multiplication on $A^F$ is associative if and only if the following compatibility axiom is satisfied: 
\beq\tag{H'4}
\sum_{(d,e,f)\in T_g}m_{de, f}(m_{d,e}(\al_d, \beta_e),\gm_f)=\sum_{(d,e,f)\in T_g}m_{d, ef}(\al_d, m_{e,f}(\beta_e,\gm_f))
\eeq
wheren $T_g$ is defined as follows. Note that $L_g$ acts diagonally on the set
\beq
\{(d,e,f)\in G\times G\times G\;|\; def=g\}
\eeq
Then by definition $T_g$ is any set of representatives for the orbits of this action of $L_g$ on the above set.
\subsubsection{}
Fix a set of representatives $\{ g_1, g_2,\cdots, g_r \}$ for the orbits of the action of $L$ on $G$. Then one can identify
\beq\lb{identific}
 A^L\simeq\oplus_{i=1}^rA(g_i)^{L(g_i)} 
 \eeq
It was noticed in \cite{scoh} that in the case the multiplication from Equation \ref{multo} is associative  then under the identification \ref{identific} this multiplication becomes:
\beq\lb{nm}
\al_i\beta_j=\sum_{x \in D}
m_{\;^{yx}g_i\;^yg_j}(c_{yx, g_i}(\al_i), c_{x, g_i}(\beta_j))
\eeq
where $y$ is chosen such that $\;^{yx}g_i\;^yg_j=g_k$ for a uniquely determined $k=k(x)$.
%%%%%%%
\subsubsection{Grothendieck rings of abelian cocentral extensions}
This type of rings were found in \cite{scoh} as the Grothendieck rings associated to abelian cocentral extensions, i.e Hopf algebras $A$ that fit the following cocentral exact sequence
\beq
k \ra k^G\ra A \ra kF \ra k
\eeq
for some finite groups $F$ and $G$.
As explained in \cite{scoh} this type of rings appear also in \cite{bouc} as the rings associated to some Green functors obtained by Dress construction from other Green rings.
\bt\lb{idr}
Let $F$ be a finite group acting coherently on a $G$-graded fusion category $\cc$. them the Grothendieck ring of $\cc^G$ is  of the same type as those introduced by Witherspoon.
\et
\bpf
Let $A=\oplus_{g \in G} A(g)$ where $A(g)=K_0(\cc^{F_g})$. Using Remark \ref{shift} one can efine $c_{x, g}: A(g)\ra A(\;^xg)$ by $[M]\mapsto [T^x(M)]$. Then it is easy to verify all the compatibility conditions from the previous subsection.
Moreover it is clear that $K_0(\cc^F):=A^F$ as defined in Equation \eqref{inf}. Then Theorem \ref{tpgr} shows that the multiplication in $K_0(\cc^F)$ coincides to the multiplication described in Equation \eqref{nm}.
\epf
\blue{ Make clear:
\bne
\item why is it enough to look at one  $x$ in each double coset in order to get all possible orbits. because otherwise one gets an element in the same orbit.
\item why my formula is equivalent to the one from corollary 2.5.
\item
why the formula from corollary is equivalent to the one from the theorem.
\ene}

\subsection{On the Grothendieck ring of the center of a fusion category} Suppose that $\cc$ is a fusion category graded by a finite group $G$
\bq
\cc=\oplus_{g \in G}\;\cc_{g}.
\eq 
Then its Drinfeld center $\cz(\cc)$ is the equivariantization of the relative center $\cz_{\cc_{1}}(\cc)$(\cite[Theorem 4.1]{gnn} by a certain action of the finite group $G$. Moreover \cite[Theorem 3.2 ]{gnn} shows that the relative center $\cz_{\cc_{1}}(\cc)$ is a $G$-crossed braided category. Applying Theorem \ref{idr} for Example \ref{gbb} it follows that the Grothendieck ring of $\cz(\cc)$ has the ring structure described above. For the Grothendieck ring of a Drinfeld double of a semisimple Hopf algebra this description was obtained in \cite{gdd}.
\bibliographystyle{amsplain}
\bibliography{bob}
\ed
\bibliographystyle{amsplain}
\bibliography{bob}
\ed
\newpage
\section{Indecomposable bimodule categories}

Suppose that $G$ acts on $\cc$ and let:
\beq
\rep(G)\ra \cc^G \ra \cc
\eeq

the associated exact sequence of categories. Let $\cc(a)$ be the full abelian subcategory generated by the simple objects induced from $(\cc_a)^{F_a}$.

Then
\bq
\cc=\opl_{a \in G}\cc(a)
\eq
is a decomposition of $\cc$ in indecomposbale $\rep(G)$ -bimodule categories.

\bp
Moreover the fusion rules are given by
\beq
\cc(a)\ot_{\rep(G)}\cc(b)=\oplus_{r}N^r_{a,b}\cc(r)
\eeq
\ep
\bpf

Note that formulae \ref{tep} implies that
\beq
S(1, M)\ot S(h, M)\simeq S(h, \res^g_H(V) \ot M)
\eeq
This implies that the functor $\cc(a)\boxtimes \cc(b) \ra \oplus_{r}N^r_{a,b}\cc(r)$ given by the tensor product from $\cc$ is $\rep(G)$-balanced. Therefore this induces a functor
\beq
\cc(a)\ot_{\rep(G)}\cc(b)\ra\oplus_{r}N^r_{a,b}\cc(r)
\eeq
which is dominant and therefore an equivalence.
\epf

\blue{study the Mueger centralizer with respect to these cosets}
\bb{\subsection{On the grothendieck group of G crossed categories}}
\newpage
\section{Actions coming from cocentral extensions of Hopf algebras}
In this section we..
\md
\subsection{Cocentral extensions of Hopf algebras}
Suppose that we have a cocentral extension of Hopf algebras 
\beq\lb{cocen}
k \ra B\xrightarrow{i} H\xrightarrow{\pi} kF\ra k
\eeq

Recall that this means the $kF^*\subset \mtc{Z}(H^*)$ via $\pi^*$. Following Proposition 3.5 of \cite{natalecoc} it follows that for any such extension one has that $F$ acts on $\Rep(B)$ and $\Rep(H)=\Rep(B)^F$.

On the other hand, using the reconstruction theorem from \cite{AD} it follows that 
\beq\lb{andc}
H \simeq B\;^{\tau}\#_{\sg} \;kF
\eeq
for some cocycle $\sg:B\ot B \ra kF$ and some dual cocyle $\tau:kF \ra B\ot B$. 
\md

Recall that in this situation the multiplication on $H$ is given by 
%\beq(a\#_{\sg}h)(b\#_{\sg}l)=a(h.b)\#hl\eeq
\bq\lb{mu}
(b\#_{\sg}f)(c\#_{\sg}g)=b(f.c)\sg(f,g)\#_{\sg}fg
\eq
and the comultiplication on $H$ is given by:
%\beq\Delta(b\#g)=(b_1\tau(g)_j\#_{\sg}g)\ot (b_2\tau(g)^j\#_{\sg}g)\eeq
\bq\lb{com}
\D(b\#_{\sg}\bar{f})=(b_1\tau(\bar{f})_i\#_{\sg}\bar{f})\ot(b_2\tau(\bar{f})^i\#_{\sg}\bar{f})
\eq
where the weak action and vocation are denoted by...
\blue{Moreover the antipode is given by the formula $$S(a\#_{\sg}g)=g^{-1}S(a)=(g^{-1}.S(a))\#_{\sg}g^{-1}.$$}
Here one may also suppose that for all $f \in F$ one has $\sg(f,f^{-1})=1$ and therefore $\bar{f}^{-1}=\bar{f^{-1}}$.  
 
Moreover, we used the same notation as in \cite{AD}, \bq\tau(f)=\tau(f)_i\ot \tau(f)^i\in B \ot B.\eq
%{\bf Check that this gives an action of $G$ on the fusion $\Rep(A)$.} 
In particular, the cocycles $\sg: F\ot F \ra U(B)$ and $\tau:F \ra B\ot B$ verify the following compatibility conditions:  

\bn{enumerate}
\item $\sg(1,g)=\sg(g,1)=1$
\item $(\eps_A\ot \id )(\tau(g))=1=(\id \ot \eps_A)(\tau(g))$
\item  \lb{cocyl} $[h.\sg(l,m)]\sg(h,lm)=\sg(h,l)\sg(hl,m)$
\item  \lb{twm} $h.(l.a)=\sg(h,l)[(hl).a]\sg^{-1}(h,l)$
\item $(\tau(g)_j)_1\tau(g)_l\ot (\tau(g)_j)_2\tau(g)^l\ot \tau(g)^j=\tau(g)_j\ot (\tau(g)^j)_1\tau(g)_r\ot  (\tau(g)^j)_2\tau(g)^r$

This means that $\tau(g)$ is a twist for $A^{op}$ or equivalently $\tau(g)^{-1}$ is a twist for $A$.
\item \lb{A} $\Delta_A(g.a) = \tau (g)(g.a_1 \ot g.a_2)\tau(g)^{-1}$ \\ (This is condition $(A)$ from \cite{AD})
\item \lb{D} $\Delta(\sg(g,h))\tau(gh)=\tau(g)(g.\tau(h)_p\ot g.\tau(h)^p)(\sg(g,h)\ot \sg(g,h))$
(This is condition $(D)$ from \cite{AD})  \end{enumerate} 
Note that conditions $(B)$ and $(C)$ from \cite{AD} are automatically satisfied.
\subsubsection{Universal gradings for Hopf algebras}
\bp
With the above notations suppose that $\cc:=\rep(B)$ and let
\begin{equation}\label{univ}
\mtc{C}=\oplus_{g \in G}\mtc{C}_g
\end{equation}
be the universal grading of $\rep(B)$ where $K(B)=kG^*$. Then $F$ acts on $G$ and the action of $F$ on $\cc$ is coherent with respect to the universal grading of $\cc$.
\ep
\bpf
First it will be shown that $F$ acts on $K(B)=kG^*$ via Hopf automorphism.
Note that if $b \in K(B)$ the $f.b \in K(B)$. Clearly $f.b \in \cz(B)$ since $F$ acts by algebra automorhisms on $B$.
\md
On the other hand using formula $(A)$ from \cite{AD} it follows that 
\beq
\Delta(f.b)=\tau(f)(f.b_1\ot f.b_2)\tau(f)^{-1}
\eeq
and if $b \in \cz(B)$ then 
\beq
\Delta(f.b)=(f.b_1\ot f.b_2).
\eeq
Thus $F$ acts by Hopf algebra automorphisms on $K(B)$. Therefore there is an action of $F$ on $G$ such that the action of $F$ on $K(B)=kG^*$ is induced by
\bq
x.p_g=p_{\;^xg}
\eq
for all $x \in F$ and $g \in G$. Now to verify  that the action of $F$ on $\rep(B)$ is coherent one has to verify that for any $M \in \co(\cc_g)$ then $T^x(M)\in \co( \cc_{\;^xg})$. Recall from Proposition 3.5 of \cite{natalecoc} that $T^x(M)=M$ as vector spaces and the action of $B$ on $T^x(M)$ is given by
\bq
b.\;^xm=(x^{-1}.b)m
\eq
Thus for any $h \in G$ one has that 
\bq
{p_h.}\;^xm=(x^{-1}.p_h)m=p_{\;^{x^{-1}}h}\;m=\delta_{\;^{x^{-1}}h, g}\;m=\delta_{h, \;^{x}g}\;m
\eq
which shows that $T^x(M)\in \co(\cc_{\;^xg})$.
\epf
\bt
Suppose that $H\simeq kF\;^{\tau}\#_{\sg} B$ is a cocentral extension of Hopf algebras as in Equation \ref{cocen}. Then $K(B)$ is a normal Hopf subalgebra of $H$.
\et
\bpf
Clearly $K(B)$ is closed under the adjoint action of $B$ on itself since it is a central Hopf subalgebra of $B$. It remains to show that $K(B)$ is closed under the adjoin action of $kF$.
Indeed one has that:
\bq
f_1bS(f_2)=(\tau(f)_i\#_{\sg}f)bf^{-1}S(\tau(f)^i)=(f.b)\tau(f)_iS(\tau(f)^i)\#_{\sg}=(f.b)
\eq  
\blue{ 
Note that
\bq
\tau(f)_iS(\tau(f)^i)=1
\eq
since $\eps(f)1=f_1S(f_2)=(\tau(f)_i \#_{\sg}f) S(\tau(f)^i\#_{\sg}f)$
}
\epf
\subsection{Clifford theory for $K(B)\subset H$}
\subsubsection{Conjugate modules}
Let $B \subset A$ be a normal Hopf subalgebra of a semisimple Hopf algebra $A$ and let $M$ be an irreducible $B$-module with associated character $\al \in C(B)$. We recall the following notion of conjugate module introduced in \cite{coset}. It was also previously considered in \cite{Schgal} in the cocommutative case.

If $W$ is an $A^*$-module then $W\ot M$ becomes a $B$-module with the following structure:
\bn{equation}\label{def}
 b(w\ot m)=w_0 \ot (S(w_1)bw_2)m
\end{equation}
for all $b \in B$, $w \in W$ and $m \in M$.
Here we used that any left $A^*$-module $W$ is a right
$A$-comodule via $T(w)=w_0\ot w_1$.
\subsubsection{Formula for conjugate module}
It can be checked that if $W \simeq W'$ as $A^*$-modules then $W\ot M \simeq W'\ot M$. Thus for any irreducible character $d \in \mtr{Irr}(A^*)$ associated to a simple $A$-comodule $W$ one can define the $B$-module $\;^dM:= W\ot M$. If $\al \in C(B)$ is the character of $M$ then the character $^d\al$ of $^dM$ is given by

\bn{equation}\label{chfom}
^d\al(x)=\al(Sd_1xd_2)
\end{equation}
for all $x \in B$ (see Proposition 5.3 of \cite{coset}).

\bl
Let $B$ be a normal Hopf subalgebra of $A$ and $M$ be an irreducible $A$-module. Then the stabilizer $Z_A(\al)$ is the largest subcoalgebra of $A$ that stabilizes $M$.
\el
\bpf
The proof follows from Proposition of \cite{doc} where it shown that the set of simple subcoalgebras stabilizing a given module is closed under product and $"*"$.
\epf
\bt
Suppose that $\tau(f)=1\ot 1$ is a trivial cocycle. Then Clifford'c correspondence holds for the above extension.
\et
\bpf
\epf
\blue{The inverse function via coset decomposition
Claims
\beq
\cz_{\cd}(\cc_g)^{F_g}=\ci_g
\eeq
\beq
T^x(M)=A^*(xF_g)\ot_{A^*(F_g)}M
\eeq
on $A^*(F_g)\btw A$-mod
Following \cite{AD} the action is
\bq
r(b_1)ar^{-1}(b_2)
\eq
}
\bpf Let $\Gm$ be a set of representatives for the left cosets of $G_Y$ inside $G$. Then one can choose $$\;^f\Gm:=\{ftf^{-1}\;| t \in \Gm \}$$
 as a set of representatives for the left cosets of $G_{\ro^f(Y)}$ inside $G$. As objects in $\cc$ one has 
\bq
 S_{Y,\pi}=\opl_{t \in \Gm}\ro^t(Y) \ot V_{\pi}  \eq and $$S_{Y, \;^f\pi}=\opl_{t\in \Gm}\ro^{ftf^{-1}}(\ro^f(Y))\ot V_{\pi}.$$ Thus one can define the isomorphism 
\bq
 \phi_f: S_{Y,\pi}\ra S_{\ro^f(Y),\;^f\pi}%\input{fusionsub-equiv3.tex}
 \eq
 which on the components is given by $$\ro^t(Y)\xra{\ro^t((c_Y^h)^{-1})} \ro^t(\ro^h(Y))\xra{\ro_2^{t,h}} \ro^{th}(Y)=\ro^{ft'}(Y)\xra{\ro_2^{ft'f^{-1},f}}\ro^{ft'f^{-1}}(\ro^f(Y))$$
 where $t'$ is chosen such that $t=ft'(\mtr{mod} G_Y)$ more precisely $th=ft'$ with $h \in G_Y$. It remains to show that $\phi_f$ is an isomorphism in $\cc^G$.
\epf
\newpage
\section{Indecomposable equivariantized objects of tensor categories}
\subsection{Semisimplicity}
\subsubsection{Group actions on Hom-spaces}
Let $M, N \in \cc^{G}$ be two equivariantized objects. There is an action of G on $\Hom_{\cc}(M, N)$ given as follows
\bq
g.f=\mu^{N}_{g}T^{g}(f)(\mu^{M}_{g})^{-1}
\eq
for all $f \in \Hom_{\cc}(M, N)$ and $g \in G$. It is easy to check that
\bq
g.(uv)=(g.u)(g.v)
\eq
for any $M, N, P\in \cc^{G}$ and any two morphisms $M\xra{v} N\xra{u} P$ in $\cc$.
\md
Note also that a morphism $f \in \Hom_{\cc}(M, N)$ is a morphism in $\cc^G$ if and only if $g.f=f$ for all $g\in G$.
\subsubsection{Canonical properties of induction and restriction functors} Let $H$ be a subgroup of $G$ and $W\in \cc^{G}$. There is an epimorphism in $\cc^{H}$
\bq
\ind^{G}_{H}(\res^{G}_{H}(W)) \xra{\al_{W}} W 
\eq
given by the  canonical projection $\oplus_{t \in G/H}T^{t}(W)\ra W$. It is easy to verify that this is a morphism in $\cc^{H}$. Indeed, for any $g \in G$ the equivariant structure of $\ind^{G}_{H}(\res^{G}_{H}(W))$ on the component $\ro^{1}(W)$ is given by:
\bq
T^{g}(T^{1}(W))=T^{g}(W)\xra{(T^{r,h}_{2})^{-1}} T^{r}(T^{h}(W))\xra{} T^{r}(W)
\eq
\md 
Moreover for any $U \in \cc^{H}$ there is also a canonical projection
\bq
\res^{G}_{H}((\ind^{G}_{H}(U))) \xra{\beta_{U}} U 
\eq
given by the canonical projection
\subsubsection{The endomorphism sum in $\mtr{End}_{\cc^{G}}(U)$}
Note that there is an endomorphism of $\ind(U)$ such that
\bq
\sum_{i}x_if=1_{U}
\eq

3. Green theorem for index p.

\subsubsection{Mashcke theorem for equivariantized tensor categories.}

\bt
The tensor category $\cc^G$ is fusion iff char k does not divide the order of $G$ and the category $\cc$ is semisimple.
\et

\bpf
Suppose that $\pi: M \ra N$ is an epim and $h$ a splitting of $\pi$ in $\cc$. Then define $h' : n \ra M$ by
\bq
h'= \frac{1}{|G|}\sum g.h
\eq

Then $h'$ is a morphism of G equivariant fusion categs since $g.h'=h'$ for all $g\ in G$. Moreover $h'$ splits $\pi$ since
\bq
\pi h'=\frac{1}{|G|}\sum_{g\in G} \pi g.h=\frac{1}{|G|}\sum_{g \in G} g(\pi h)=\frac{1}{|G|}
\sum_{g \in G} g. 1_{N}=1_{N}
\eq

The converse also follows. $F: \cc^G \ra \cc$ is a surjective tensor functor. Then $\cc$ is semisimple. Moroever $\rep(G)$ is a tensor subcategory of $\cc$ and thus semisimple. Then Mashke's theorem implies that  order of $G$ does does not divide the characteristic of $k$.
\epf
\blue{
There is an underlined category $F(\cc^g)$  which is semisimple.}
\md
5. Indecomposable modules for equiv tensor categories.
\md
6. Define blocks and do evth as for modular categories.
\md
7. Compact groups or of finote index both of them.

In this section we study

\subsection{Vertices and sources}

\subsection{Relatively $H$-free objects}
 is relatively $H$-free if for any object 
\subsection{Relatively $H$-projective objects}
%\bn{def}
Let $(W, \mu)\in \cc^{G}$ and $H$ a subgroup of $G$. Then $W$ is
called $H$-projective if, whenever there is an exact sequence of the form 
\beq
0\ra V\ra U\ra W\ra 0
\eeq
such that under restriction to $H$ the exact sequence splits then it also splits in $\cc^{H}$.
%\end{def}
\md
\bl For any two $B$-modules $P_1$ and $P_2$ the direct sum $P :=P_1 \oplus P_2$ is relatively projective if and only if $P_1$ and $P_2$ both are relatively projective.
\el
We generalize the following characterization of $H$-projectivity due to Higman for representations of groups.
\bt Let H be a subgroup of the group G. Then each of the following properties for $ W\in \cc^{G}$ are equivalent:
\md
(a) $W$ is H-projective;
\md
(b) $W|\ind^{G}_{H}(W)$
\md
(c) $W|\ind^{G}_{H}(U)$ for some object $U\in \cc^{G}$ 
\md
(d) there is $f\in \mtr{End}_{\cc^{H}}(W)$ such that $\sum_{i}x_{i}.f=1_{W}$.
\et
\bpf
$a)\implies b)$ There is a surjective morphism $\ind^{G}_{H}(W)\ra W$ define as... that splits when restricted to $\cc^{H}$.
\md Trivially $b)\implies c)$.
\md
$c)\implies a)$
\epf
%(d) there exists f~ End,,( W) such that C;= , x; ' f x i = 1,
if possible to write Mackey for non semisimple fusion categories

\subsection{Vertices and sources}
Let $V$ be an indecomposable $G$-equivariant object. Then a subgroup $Q$ of $G$ is a
vertex of $V$ if $V$ is $Q$-projective but $V$ is not $H$-projective for any proper subgroup of $Q$. 
\md

\bl
Let p be a prime number and suppose that all prime
and let  $k$ be a  field of characteristic
p ; then all vertices are p-groups.
\el
\bpf Let $H_1 \in V_0(M)$ for some R[G]-module M and let $H_0 \subset H_1$ be a
p-Sylow subgroup of $H_1$. We claim that M is relatively $k[H_0]$-projective
which, by the minimality of $H_1$, implies that $H_0 = H_1$ is a $p$-group. Let
Let $L\ra N$ be an epimorphism that splits in 
$\cc^{H_0}$ 
\md
Since $[H_1 : H_0]$ is invertible
in R it follows from Prop. 18.8 that $M$ viewed as an $R[H_1]$-module
is relatively $R[H_0]$-projective. Hence there exists an $R[H1]$-module homomorphism $\al_1:M \ra L$ such that $\al_1$ is a splitting for $\beta$. But $H_1 \in V(M)$, so there
further must exist an $R[G]$-module homomorphism $\al: M \ra L$ satisfying $\beta \al=\id$.
\epf
It follows from Theorem 5.l(c) that corresponding to any vertex
$Q$ of $V$ there is at least one indecomposable $Q$-equivariant object $U$ such that $V| U^{G}$ ;
such an Q-equivariant object is called a source of $V$.
\subsection{Green correspondence}

Throughout this subsection we assume the commutative ring $R$ to be noetherian
and complete and $R= Jac(R)$ to be artinian. 
\md
We fix a subgroup $H\subset G$.
We consider any finitely generated indecomposable $R[G]$-module $M$ such
that $H \in V(M)$. Let
\bq
\res^g_H(M) = L_1\oplus\cdots\oplus L_r
\eq
be a decomposition of $M$ as an $R[H]$-module into indecomposable $R[H]$-
modules $L_i$. 
%%%
\bc Let $H$ and $L$ be subgroups of the group G.
(i) If $x^{{-1}}Hx\subset  L$ for some $x \in G$, then each $H$-projective $G$-equivariant object is $L$-projective.
(ii) IF $H \subset $L and $V$ is an H-projective $L$-equivariant object, and W is an AC-equivariant object such that $W | V^{G}$, then $W$ is an $H$-projective $G$-equivariant object.
\ec
\bpf(i) Let $W$ be an $H$-projective $G$-equivariant object. By part (c) of the
theorem there exists an AH-equivariant object U such that W I UG.T hen U �3 x is an
A(x- 'Hx)-subequivariant object of U" and it is easily seen that (U �3 x)" = UG as
$G$-equivariant objects. Let V be the induced AL equivariant object (U @I x)". Since inducing is a
transitive operation, p = (U �3 x)" = UG and so W I p. Thus the theorem
shows that W is Lprojective.
(ii) Since V is an H-projective AL-equivariant object, part (c) of the theorem shows
that $V | U^{G}$ for some AH-equivariant object $U$. Then $V|  ( U^{L})^{G} = U^{G}$, so $W | U^{G} $ and hence $W$ is $H$-projective.
\epf
\md
Let H be a normal subgroup of index $p$ in the group $G$, and
suppose that k is a perfect field of characteristic p. Let W be an absolutely
indecomposable kH-module, and suppose that the indecomposable components
of (Wc)H are all isomorphic to W Then Wc is an absolutely
indecomposable kG-module.

\bibliographystyle{amsplain}
\bibliography{bob}
\ed

One has that 
\beqarn
S({g, M})\ot S({h,N}) & = & (\oplus_{r \in F/F_g}T^r(M))\ot (\opl_{s \in F/F_h}T^s(N))\\  & \simeq & \opl_{r \in F/F_g, s\in F/F_h}T^r(M)\ot T^s(N)
\eeqarn

{\center
\begin{tikzpicture}[>=latex]\lb{hexagon}
    \def\radius{5cm} % change to an appropriate value
    \node (h0A) at (50:\radius)   {$T^{xhx\inv}(T^{xlx\inv}(T^{x}(V))$};
    \node (h0C) at (20:\radius)    {$T^{xhx\inv}(T^{xl}(V))$};
    \node (h1B) at (-30:\radius)  {$M\ot (X\ot Y)$};
    \node (h1A) at (-150:\radius) {$(X\ot Y)\ot M$};
    \node (h1C) at (180:\radius)  {$X\ot (Y\ot M)$};
    \node (h0B) at (150:\radius)  {$X\ot (M \ot Y)$};
    \node {l0A} at (150:\radius) {Me}
    \path[->,font=\small]
        (h0A) edge node[auto] {$\gm_X\ot 1$} (h0C)
        (h0C) edge node[auto] {$\al_{M,X,Y}$} (h1B)
        (h1B) edge node[auto] {$\gm_{X\ot Y}$} (h1A)
        (h1C) edge node[auto] {$\al^{-1}_{X,Y,M}$} (h1A)
        (h1C) edge node[auto] {$1\ot \gm_Y$} (h0B)
        (h0B) edge node[auto] {$\al^{-1}_{X,M,Y}$} (h0A);
\end{tikzpicture}

{\SMALL
\begin{tikzpicture}
  \matrix (m) [matrix of math nodes, row sep=4em,column sep=8em,minimum width=5em]
  { T^{xax\inv}(T^{x}(M\ot N)) &     T^{xax\inv}(T^{x}(M)\ot T^{x}(N))\\%1
  T^{xa}(M\ot N) &  T^{xax\inv}(T^{x}(M))\ot T^{xax\inv}(T^{x}(N))\\%2
  T^{x}(T^{a}(M\ot N)) & T^{xa}(M)\ot T^{xa}(N)\\%3
  T^{x}(T^{a}(M)\ot T^{a}(N)) & T^{x}(T^{a}(M))\ot T^{x}(T^{a}(N))  \\%4
  T^{x}(M\ot N) & T^{x}(M)\ot T^{x}(N)\\};%5
  \path[-stealth]
  (m-1-1) edge node [above] {$T^{xax\inv}((T^{x}_{2})^{M, N})$} (m-1-2)
    (m-5-1) edge node [above] {${(T^{x}_{2})}^{M, N}$} (m-5-2)
         
         (m-2-1) edge node [left] {$(T^{x,a}_{2})^{-1}_{M\ot N}$} (m-3-1)
         (m-3-1) edge node [left] {$T^{x}((T_{2}^{a})^{M, N})$} (m-4-1)
         (m-4-1) edge node [left] {$T^{x}(\mu^{M}_{a}\ot \mu^{N}_{a})$} (m-5-1)
         
            (m-2-2) edge node [right] {$(T_{2}^{xax\inv, x})_{M}\ot (T_{2}^{xax\inv, x})_{N}$} (m-3-2)
         (m-3-2) edge node [right] {$(T_{2}^{x, a})^{{-1}}_{M}\ot (T_{2}^{x, a})^{{-1}}_{N}$} (m-4-2)
         (m-4-2) edge node [right] {$ T^{x}(\mu^{M}_{a}) \ot T^{x}(\mu^{N}_{a})$} (m-5-2)
               (m-1-1) edge node [left] {$(T^{xax\inv, x}_{2})_{M\ot N}$} (m-2-1)
    (m-1-2) edge node [right] {$(T^{xax\inv}_{2})^{{T^{x}(M), T^{x}(N)}}$} (m-2-2)
    (m-2-1) edge [dashed] node [above] {$(T^{xa}_{2})^{M,N}$} (m-3-2)
     (m-4-1) edge [dashed] node [above] {$({T_{2}^{x})}^{{T^{a}(M)},{T^{a}(N)}}$} (m-4-2)
    ;
\end{tikzpicture}
}